\documentclass[10pt,twoside]{article}
\usepackage[english,russian]{babel}
\usepackage{droid}
\usepackage{famems}
\usepackage{braket}
\usepackage{cmap}
\usepackage{graphics}
\usepackage{graphicx}
\usepackage{graphpap}
\usepackage{amsmath,amstext,amsfonts,amssymb,amsthm}
\usepackage{mathtools}
\usepackage{dsfont}
\usepackage[mathscr]{eucal}
\usepackage{array}
\usepackage{ctable}
\usepackage{multicol}
\usepackage{multirow}
\usepackage{multicol}
\usepackage{tabularx}
\usepackage{xcolor}
\usepackage{color}
\usepackage{colortbl}
\usepackage{bibunits}
\usepackage{ifthen}
\usepackage{epstopdf}
\usepackage[normalem]{ulem}
\usepackage{rotating}
\usepackage{hhline}
\usepackage{pdflscape}
\usepackage{fontawesome}
\usepackage{scrextend}
\usepackage{arydshln}
\usepackage{tikz}
\usepackage{pict2e}
\usepackage{url}
\usepackage[breaklinks]{hyperref}
\usepackage{hyperref}
\hypersetup{
    colorlinks,
    citecolor=black,
    filecolor=black,
    linkcolor=black,
    urlcolor=
teal
}

\makeatletter
\renewcommand{\paragraph}{\@startsection{paragraph}{4}{0ex}%
   {-3.25ex plus -1ex minus -0.2ex}%
   {1.5ex plus 0.2ex}%
   {\normalfont\normalsize\tt}}
\makeatother
\begin{document}
\baselineskip=13pt

%%%%%%%%%%%%%%%%%%%%%%%%%%%%%%%%%%%%%%%%%%%%%%%%%%%%%%%%%%%%%%%%%%%%%%%%%%%%%%%%%%%%%%%%%%%%%%%%%%%%%%%%%%%%%%%%%%%%%%%%%%%%%%%%%%%%%%%%%%%%%%%
%%%%%%%%%%%%%%%%%%%%%%%%%%%%%%%%%%%%%%%%%%%%%%%%%%%%%%%%%%%%%%%%%%%%%%%%%%%%%%%%%%%%%%% vorobyev %%%%%%%%%%%%%%%%%%%%%%%%%%%%%%%%%%%%%%%%%%%%%%
%%%%%%%%%%%%%%%%%%%%%%%%%%%%%%%%%%%%%%%%%%%%%%%%%%%%%%%%%%%%%%%%%%%%%%%%%%%%%%%%%%%%%%%%%%%%%%%%%%%%%%%%%%%%%%%%%%%%%%%%%%%%%%%%%%%%%%%%%%%%%%%
\newcounter{ctrwar}\setcounter{ctrwar}{0} % 12 Jan 2017
\newcounter{ctrdef}\setcounter{ctrdef}{0}
\newcounter{ctrdefpre}\setcounter{ctrdefpre}{0}
\newcounter{ctrTh}\setcounter{ctrTh}{0}
\newcounter{ctrnot}\setcounter{ctrnot}{0}
\newcounter{ctrATT}\setcounter{ctrATT}{0} % 12 Jan 2017
\newcounter{ctrcor}\setcounter{ctrcor}{0}
\newcounter{ctrAx}\setcounter{ctrAx}{-1}
\newcounter{ctrexa}\setcounter{ctrexa}{0}
\newcounter{ctrlem}\setcounter{ctrlem}{0}
\newcounter{ctrPRO}\setcounter{ctrPRO}{0}
\newcounter{ctrrem}\setcounter{ctrrem}{0}
\newcounter{ctrass}\setcounter{ctrass}{0}
\newcounter{ctrlist}\setcounter{ctrlist}{0}

\newcommand{\bfPB}{\mbox{\protect\reflectbox{\bf P}\hspace{-0.4em}{\bf B}}}
\newcommand{\bfEl}{\mbox{\protect\reflectbox{\bf E}}}
\newcommand{\bfEr}{\mathbf{E}}
\newcommand{\bfEE}{\mbox{\protect\reflectbox{\bf E}\hspace{-0.4em}{\bf E}}}
\newcommand{\bfphi}{\mbox{\boldmath$\varphi$}}
\newcommand{\bfPhi}{\mbox{\boldmath$\Phi$}}
%%%%%%%%%%%%%%%%%%%%%%%%%%%%%%%%%%%%%%%%%%%%%%%%%%%%%%%%%%%%%%%%%%%%%%%%%%%%%%%%%%%%%%%%%%%%%%%%%%%%%%%%%%%%%%%%%%%%%%%%%%%%%%%%%%%%%%%%%%%%%%%
%%%%%%%%%%%%%%%%%%%%%%%%%%%%%%%%%%%%%%%%%%%%%%%%%%%%%%%%%%%%%%%%%%%%%%%%%%%%%%%%%%%%%%% vorobyev %%%%%%%%%%%%%%%%%%%%%%%%%%%%%%%%%%%%%%%%%%%%%%
%%%%%%%%%%%%%%%%%%%%%%%%%%%%%%%%%%%%%%%%%%%%%%%%%%%%%%%%%%%%%%%%%%%%%%%%%%%%%%%%%%%%%%%%%%%%%%%%%%%%%%%%%%%%%%%%%%%%%%%%%%%%%%%%%%%%%%%%%%%%%%%

\numberwithin{equation}{section}

\renewcommand{\figurename}{\scriptsize  Figure}
\renewcommand{\tablename}{\scriptsize  Table}
\renewcommand\refname{\scriptsize References}
\renewcommand\contentsname{\scriptsize Contents}
\setcounter{page}{25}

\titleOneAuthorenOne{
Postulating the theory of experience and  chance\\ as a theory of co$\sim$events (co$\sim$beings)

 }{Oleg Yu. Vorobyev}{
Institute of mathematics and computer science\\
Siberian Federal University \\
Krasnoyarsk \\
\tiny
\url{mailto:oleg.yu.vorobyev@gmail.com}\\
\url{http://www.sfu-kras.academia.edu/OlegVorobyev}\\
\url{http://olegvorobyev.academia.edu}}

\label{vorobyev2}

\colontitleen{Oleg Yu Vorobyev. Postulating the theory of experience and  chance} \footavten{O.Yu.Vorobyev} \setcounter{footnote}{0}
\setcounter{equation}{0} \setcounter{figure}{0} \setcounter{table}{0} \setcounter{section}{0}

\vspace{-22pt}

\begin{abstracten}
\textit{The aim of the paper is the axiomatic justification of the theory of experience and  chance, one of the dual halves of which is the
Kolmogorov probability theory. The author's main idea was the natural inclusion of Kolmogorov's axiomatics of probability theory in a number of
general concepts of the theory of experience and  chance. The analogy between the measure of a set and the probability of an event has become
clear for a long time. This analogy also allows further evolution: the measure of a set is completely analogous to the believability of an event.
In order to postulate the theory of experience and  chance on the basis of this analogy, you just need to add to the Kolmogorov probability
theory its dual reflection --- the believability theory, so that the theory of experience and  chance could be postulated as the certainty
(believability-probability)  theory on the Cartesian product of the probability and believability spaces, and the central concept of the theory is
the new notion of co$\sim$event as a measurable binary relation on the Cartesian product of sets of elementary incomes and elementary outcomes.
Attempts to build the foundations of the theory of experience and  chance from this general point of view are unknown to me, and the whole range
of ideas presented here has not yet acquired popularity even in a narrow circle of specialists; in addition, there was still no complete system of
the postulates of the theory of experience and  chance free from unnecessary complications. Postulating the theory of experience and  chance
can be carried out in different ways, both in the choice of axioms, and in the choice of basic concepts and relations. If one tries to achieve the
possible simplicity of both the system of axioms and the theory constructed from it, then it is hardly possible to suggest anything other than
axiomatization of concepts co$\sim$event and its certainty (believability-probability). The main result of this work is the \textbf{axiom of
co$\sim$event}, intended for the sake of constructing a theory formed by dual theories of believabilities and probabilities, each of which itself
is postulated by its own Kolmogorov system of axioms. Of course, other systems of postulating the theory of experience and  chance can be
imagined, however, in this work a preference is given to a system of postulates that is able to describe in the most simple manner the results of
what I call an experienced-random experiment.}
\end{abstracten}\\[-21pt]
\begin{keywordsen}
\emph{Eventology, event, co$\sim$event, experience, chance, to experience, to happen, to occur, theory of experience and  chance, theory of
co$\sim$events, axiom of co$\sim$event, probability, believability, certainty (believability-probability), probability theory, believability
theory, certainty theory.}
\end{keywordsen}

\clearpage
\section{Co$\sim$event as a set of dual pairs}%\\  \tiny$\sim$ \textsf{0000/vorobyev2-0000, 4 Feb 2017}}

{\footnotesize
\hfill \emph{The solution to any problem}\\[-16pt]

\hfill \emph{begins with correcting the names.}\\[-16pt]

\hfill \emph{If the names are wrong,}\\[-16pt]

\hfill \emph{then the concepts have no basis.}\\[-16pt]

\hfill \emph{If the concepts have no basis,}\\[-16pt]

\hfill \emph{then the events can not occur.}

}

%\smallskip

{\scriptsize
\hfill Confucius (551 -- 479 BC)\\[-14pt]

\hfill Analects (Edited Conversations)\\[-14pt]

\hfill Chapter XIII, 3.

}

%\bigskip

I'll start with one detail at which you should linger. Among the reasons that gave rise to the theory of experience and chance, for a long time it would be possible to linger on the philosophy of the duality of being. But our milestone is quite different.

The ``tacit'' Kolmogorov axiom defines each event $x$ as some subset $x \subseteq \Omega$ of elementary outcomes $\omega \in
\Omega$ such that
\begin{itemize}
\item when a one $\omega \in x$ happens,\\[-1pt] they say that the \emph{\textbf{event $x$ happens}};
\item otherwise, when no elementary outcome $\omega \in x$ happens,\\[-1pt]
      they say that the \emph{\textbf{event $x$ does not happen}}.
\end{itemize}
So, the fact that the \emph{event $x$ happens for an elementary outcome $\omega\in\Omega$} is defined by the ``tacit'' Kolmogorov axiom as a realization of the membership relation: $\omega\in x$ (see Axiom \ref{Ax-event-Kolmogorov} on page \pageref{Ax-event-Kolmogorov}). For reasons unknown to us, this postulate is not included in Kolmogorov's axiomatics of probability theory explicitly: it received from its creator the role of only a preliminary definition.

At the same time, it is this statement, as the \emph{axiom of the event}, that can serve as an essential aid in delimiting probability theory and general measure theory. Moreover, in the new \emph{theory of experience and chance (TEC)} this \emph{axiom of the event} enters as one of the dual halves in the \emph{axiom of the co$\sim$being} (see Axiom \ref{Ax-co-event} on page \pageref{Ax-co-event}), without explicit support for which the new theory can not take place because the TEC sees in everything, that we have always understood under the events, dual pairs\footnote{
This one and a number of subsequent formulas use the bra-ket terminology and bra-ket notation $\braket{\cdot|\cdot}$, which are defined below and which largely rely on what I call an \emph{element-set labelling} (see \cite{Vorobyev2016famems1}).}:
\begin{equation}\label{def-dual-pair}
\mbox{$\Big\langle\mbox{bra-event (experience)}\,\Big|\,\mbox{ket-event (chance)}\Big\rangle$.}
\end{equation}
and defines its central concept, \emph{co$\sim$event (experience$\sim$chance)}, as the set of such dual pairs.

The definition of co$\sim$event as a set of dual pairs (\ref{def-dual-pair}) is not someone's whim, and certainly not mine. I will venture to say that this is only the ``wish'' of Kolmogorov's theory of probability, which despite its perseverance is still hidden from prying eyes. And the point is this.

It suffices to imagine a finite set of Kolmogorov events $\frak X\subset\mathcal A$, chosen from the sigma-algebra of the probability space  $(\Omega,\mathcal{A},\mathbf{P})$, which consists of Kolmogorov events $x \in \frak X$, defined in Kolmogorov's theory of probability as measurable subsets $x \subseteq \Omega$ of elementary outcomes  $\omega \in \Omega$; so that before your eyes there is such the following chain of two relations of membership:
\begin{equation}\label{membership-inclusion}
\omega \in x \in \frak X ,
\end{equation}
where the Kolmogorov event $x \subseteq \Omega$ acts in dual roles: an \textbf{element of the set} $\frak X$, and a \textbf{subset of the set} $\Omega$.

Such a dilemma is not only uncommon in a hard-to-see corpus of mathematical theories using the language of set theory, but rather it is truism. But in probability theory, this truism has proved to be a natural carrier of the deep sense of definition (\ref{def-dual-pair}), which suggests working with each concept of the Kolmogorov theory of events and their probabilities as with a dual pair consisting of an experience of observers and an observation of chance.
As a result of such \emph{element-set duality}, the Cartesian product
\begin{equation}\label{dual-pair}
\braket{\Omega|\Omega} = \bra{\Omega}\times\ket{\Omega},
\end{equation}
of the \emph{bra-set $\bra{\Omega}$ (the set of experiences of observes)} and the \emph{ket-set $\ket{\Omega}$ (the set of chances of observation)} becomes the mathematical model of $\Omega$, and the dual pair
\begin{equation}\label{dual-pair}
\braket{x|x}\subseteq\braket{\Omega|\Omega},
\end{equation}
becomes the mathematical model of each \emph{event} $x\subseteq\Omega$ as the \emph{co$\sim$event}.
The first element of the pair, the bra-event $\bra{x}\subseteq\bra{\Omega}$, plays a role of the event $x$ as an \textbf{element of the set} $\frak X$ and describes an \emph{experience of observer of what $x$ was}, and the second element, the ket-event $\ket{x}\subseteq\ket{\Omega}$, plays a dual role of the event $x$ as a \textbf{subset of the set} $\Omega$ and describes an \emph{observation of what $x$ is}.

Moreover, a \emph{duality of an element and a subset} \cite{Vorobyev2016famems1}, which naturally manifests itself in the concept of the Kolmogorov event in probability theory, ensures the continuation of the chain of two membership relations (\ref{membership-inclusion}) to the following:
\begin{equation}\label{membership-inclusions-3}
\omega \in x \in X \in \reflectbox{\bf S}^{\frak X}\subseteq\mathscr{P}(\frak X),
\end{equation}
where now a subset of the Kolmogorov events $X \subseteq \frak X$ also appears in the dual role as an \textbf{element of the set} $\reflectbox{\bf S}^{\frak X}\subseteq\mathscr{P}(\frak X-\{\varnothing^\Omega\})$, and a \textbf{subset of the set} $\frak X$. This time another dual pair:
\begin{equation}\label{dual-terraced-pair}
\Braket{\textsf{Ter}_{X/\!\!/\frak X}|\textsf{ter}(X/\!\!/\frak X)}\subseteq\braket{\Omega|\Omega}
\end{equation}
becomes the mathematical model of each so called  \emph{terraced event numbered by $X$} as the \emph{terraced co$\sim$event}.
The left element of the pair, \emph{terraced bra-event}
\begin{equation}\label{terraced-bra}
\Bra{\textsf{Ter}_{X/\!\!/\frak X}}=\sum_{x\in X} \bra{x}\subseteq\bra{\Omega},
\end{equation}
numbered by $X\subseteq\frak X$, \textbf{as by subset of the set} $\frak X$, is defined by the \emph{union of subset of experiences of observes} $\bra{x}, x \in X$, and the right element, \emph{terraced ket-event}
\begin{equation}\label{terraced-ket}
\ket{\textsf{ter}(X/\!\!/\frak X)}=\bigcap_{x\in X} \ket{x} \bigcap_{x\in \frak X -X} \ket{x}^c\subseteq\ket{\Omega},
\end{equation}
numbered by $X\in\reflectbox{\bf S}^{\frak X}$, \textbf{as by element of the set} $\reflectbox{\bf S}^{\frak X}$, is defined as
the \emph{observation of intersection of the set of chances} $\ket{x}, x \in X$ and $\ket{x}^c=\ket{\Omega}-\ket{x}, x \in \frak X -X$, where
$\ket{x}^c=\ket{\Omega}-\ket{x}$ is a complement of the ket-event $\ket{x}$ to the ket-set $\ket{\Omega}$.

Although the previous preliminary text ``slightly'' runs ahead and contains some mathematical misunderstandings due to the premature use of the still-unknown bra-ket concepts and notations of the element-set labelling, but we will still have time and the possibility of their correct definition to show convenience, practicality and unbearable fruitfulness of dual mathematical models co$\sim$event (\ref{dual-pair}) and terraced co$\sim$event (\ref{dual-terraced-pair}) as dual pairs that are unusually effective not only in theory but also in applications.

\section{Warnings}%\\  \tiny$\sim$ \textsf{0000/vorobyev2-0000, 19 March 2017}}

\texttt{Warning \!\refstepcounter{ctrwar}\arabic{ctrwar}\,\label{war-membership1}\itshape\footnotesize (dual interpretation of a chain of memberships)\!.} I consider it my duty to warn the reader of a perfectly understandable desire not to detain a glance at the chain of three relations of memberships
\begin{equation}\label{membership-inclusions-3a}
\omega \in x \in X \in  \reflectbox{\bf S}^{\frak X} \subseteq \mathscr{P}(\frak X),
\end{equation}
which for all and at once seems to be not worthy of attention as any other set-theoretical banality (see Warning \ref{war-membership2}).
In fact, (\ref{membership-inclusions-3a}) does not contain any set-theoretic news.
However, one news for the theory of experience and chance and, in particular, for the probability theory in it all the same is:
the dual interpretation of the chain of membership relations rightfully plays, in my opinion, a key role in defining the basic concepts of TEC.
This role is so key that this duality is considered by me as the basis of the axiomatics of the new theory. The ordinary chain of membership relations, a simple sequence of binary relations of elements, subsets and sets of subsets, serves as an inevitable set-theoretic cause that forces the modern understanding of Kolmogorov's theory of probability to be transformed from an important but special case of general measure theory to one of the dual halves of the new theory of experience and chance.
This news is naive as a too literal adherence to the meaning of binary relations between set-theoretic concepts, and is revealing, revealing something that is still unknown, in the theoretical disclosure of our own practices of observers and observations. And sometimes, as in the study of experience and chance, such disclosure turns out to be an incomprehensible exact theory.

\texttt{Warning \!\refstepcounter{ctrwar}\arabic{ctrwar}\,\label{war-membership2}\itshape\footnotesize (membership relations and paradoxes of naive set theory)\!.} \, Some mathematical relations such as ``member of'' and ``subset of'', generally speaking, should not be understood as  binary relations because its domains and codomains cannot be sets in usual systems of axiomatic set theory. For example, if you
try to model the general concept of membership as a binary relation ``$\in$'', then
then for this you will have to define the domain and the codomain, which can be a class of all sets.
But such a class is not a set in the naive set theory, and the assumption that the relation ``$\in$'' is defined on all sets leads to a contradiction from the well-known Russell paradox. At the same time, in the overwhelming majority of mathematical contexts, links to the relation ``member of'' and ``subset of'' are absolutely harmless, because they are tacitly limited to some set which is clear from the context.
The removal of this problem consists in choosing each time a sufficiently large set $A$, which contains all objects of interest, and work with the restriction ``$\in_A$'' instead of ``$\in$''.
Similarly, the relation ``$\subseteq$'' must also be limited to the relation ``$\subseteq_A$'' to have some domain  $A$ and the codomain $\mathscr {P}(A)$, set of all subsets of $A$. Therefore, the chain of three membership relations (\ref{membership-inclusions-3a}) will always be understood by me as
\begin{equation}\label{membership-inclusions-3b}
\omega \ \in_\Omega \ x \ \in_{\frak X} \ X \ \in_{\reflectbox{\scriptsize\bf S}^{\frak X}} \ \reflectbox{\bf S}^{\frak X} \  \subseteq_{\mathscr{P}(\frak X)} \ \mathscr{P}(\frak X),
\end{equation}
the chain of \emph{limited by default} membership relations.

\texttt{Warning \!\refstepcounter{ctrwar}\arabic{ctrwar}\,\label{war-membership3}\itshape\footnotesize (relative subsets and relative empty subsets)\!.} Since in the theory of experience and chance one has to deal simultaneously with subsets of sets of different levels, we will need unusual, but convenient notation, directly indicating what subsets of which set is spoken.
For example, if we are talking about subsets $x\subseteq \Omega$, $X\subseteq\frak X$, or $\mathscr{O}\subseteq\mathscr{P}(\frak X)$, then denotations of subsets $x,X$, or $\mathscr{O}$, when appropriate, we will write more fully: $x/\!\!/\Omega$, $X/\!\!/\frak X$, or $\mathscr{O}/\!\!/\mathscr{P}(\frak X)$, directly specifying in which sets these subsets contain.
Especially we will have to deal with empty subsets: $\varnothing/\!\!/\Omega$, $\emptyset/\!\!/\frak X$, or $\emptyset/\!\!/\mathscr{P}(\frak X)$, for which we introduce more compact notation: $\varnothing^\Omega=\varnothing/\!\!/\Omega$, $\emptyset^{\frak X}=\emptyset/\!\!/\frak X$, or $\emptyset^{\mathscr{P}(\frak X)}=\emptyset/\!\!/\mathscr{P}(\frak X)$, we will talk about them as \emph{relatively empty subsets}, and call \emph{$\Omega$-empty, $\frak X$-empty, or $\mathscr{P}(\frak X)$-empty subsets} correspondingly.

\texttt{Warning \!\refstepcounter{ctrwar}\arabic{ctrwar}\,\label{war-happen-experience-occur}\itshape\footnotesize (to happen, to be experienced, to occur)\!.} Theory of experience and chance, or theory of certainties, is a theory of co$\sim$events. It is a synergy of two interrelated dual theories --- the theory of believabilities and the theory of probabilities that study two dual faces of the co$\sim$event --- a ket-event, which can \textbf{\emph{happen or not happen}}, and a bra-events, which can \textbf{\emph{be experienced or not be experienced}}, in order to
the co$\sim$event itself could \textbf{\emph{occur or not occur}}.
For a long time I selected the words \emph{to happen, to be experienced, to occur} to describe the way of existence of a co$\sim$event and its dual faces.
It is possible that my choice to someone seems not entirely successful.
However, these words, in my opinion, are most similar to expressing two dual parts of what could previously be expressed in one word: \emph{to occur}.
In the theory of co$\sim$events the expression \emph{``to occur''} is understood as \emph{``to be experienced what happens''} and is associated only with a co$\sim$event, and for its dual parts ``new'' terms: \emph{to happen} for ket-events, and \emph{to be experienced} for bra-events, are fixed.
I could not find these three words right away, which helped me in the selection process to make myself forget and to ask the reader now to try to forget that the words \emph{to happen, to be experienced, to happen} are usually perceived, rather, as synonyms for each other.
This is important because in this text I intend to use them exclusively as three different mathematical terms, denoting three different concepts. Of course, this will make the style of the presentation much more difficult, but I'm ready to sacrifice the style for the sake of accuracy of expressing the main idea of the new theory about dual nature of co$\sim$event: \emph{``something occurs when one is experienced what happens''} (See Axiom \ref{Ax-co-event} on page \pageref{Ax-co-event}).

\section{``Element-set coordinates'' generated by a binary relation}%\\  \tiny$\sim$ \textsf{0000/vorobyev2-0000, 19 March 2017}}

Our goal is to divide each concept of the theory of experience and chance into two dual parts and present it in the form of a conveniently written dual pair. For the recording of such dual pairs, we are proposing, for the time being, only formally to borrow the Dirac notation \cite{Dirac1939, Dirac1964}, which are quite suitable for our purposes and well-proven in quantum mechanics. In order to continue the study of the duality of elements and sets in bra-ket notations, it is necessary to begin with the definition of some preliminary terminological set-theoretic constructions necessary for constructing the bra-ket presentation of the new theory. It is a question of the notion of a measurable binary relation as the most suitable applicant for the mathematical model of an event as a dual pair. It turned out that the measurable binary relation has very convenient labelling properties \cite{Vorobyev2016famems1}. The point is that for work in a set-theoretic space whose objects of interest serve simultaneously \emph{space} elements, \emph{sets of elements}, and \emph{sets of subsets of elements}, it is necessary to have in stock a certain coordinate system suitable for labelling both the space itself and its parts. Here, in my opinion, a slightly peculiar but effective system of set-theoretic coordinates, generated by the measurable binary relation and quite based on some labelling set $\frak X$ and some set  $\reflectbox{\bf S}^{\frak X}\subseteq\mathscr{P}(\frak X)$ of its labelling subsets, and also on the \emph{M-complement}\footnote{\label{M-complement-set}The set $\frak X ^{\!(\!c\!)}=\{x^c\colon x\in\frak X \}$ is called a \emph{complement by Minkowski (an M-complement)} of the set $\frak X$.} $\frak X ^{\!(\!c\!)}$ of the labelling set $\frak X$ and on the one-to-one corresponding $\reflectbox{\bf S}^{\frak X}$ set of its labelling subsets  $\reflectbox{\bf S}^{\frak X ^{\!(\!c\!)}}=\left\{X^{\!c(\!c\!)}\colon X\in\reflectbox{\bf S}^{\frak X}\right\}\subseteq\mathscr{P}\left(\frak X ^{\!(\!c\!)}\right)$.

Consider the \emph{measurable space} $(\Omega,\mathcal{A})$ composed of some set $\Omega$ and a sigma-algebra $\mathcal{A}$ of its
subsets and we emphasize that: \emph{elements} $\omega \in \Omega$; \emph{measurable subsets} $x \subseteq \Omega$; \emph{some set}
$\frak X =\{x\colon x\in\mathcal{A} \} \subseteq \mathcal{A}$, composed from measurable subsets $x\in\frak X$; and \emph{some set
$\reflectbox{\bf S}^{\frak X}\subseteq\mathscr{P}(\frak X-\{\varnothing^\Omega\})$  of subsets}  $X\subseteq\frak X$, consisting from measurable subsets  $x\in
X\subseteq\frak X$; \emph{until they have no meaningful interpretation} and form only a basis {\boldmath$\Lambda$} peculiar
element-set \emph{\textbf{labels}} $\lambda \in$ {\boldmath$\Lambda$} (\emph{tags, dockets, tickets}, or \emph{names}),
intended for a element-set labelling, or a nominating the parts and details of the construction that we are going to propose in the theory of experience and chance as a mathematical model of an event as a dual pair.

\texttt{\linebreak\indent Predefinition \!\refstepcounter{ctrdefpre}\arabic{ctrdefpre}\,\label{defpre-label}\itshape\footnotesize (Basic element-set labels)\!.} \emph{Basic element-set labels} $\lambda\in\mbox{\boldmath$\Lambda$}$ are called as elements,
sets and sets of subsets of the measurable space $(\Omega,\mathcal{A})$, and also results of \emph{terraced set-theoretic operations} over them, equipped with their own titles.

We'll fill up the stock of {\boldmath$\Lambda$} tags with one more label, Cartesian product
\begin{equation}\label{0000-XsX-cartesian}
\frak X  \times \reflectbox{\bf S}^{\frak X} = \left\{ (x,X) : x \in \frak X , X \in \reflectbox{S}^{\frak X} \right\},
\end{equation}
which defines a binary relation
\begin{equation}\label{0000-XsX-relation}
\mathscr{R}_{\frak X ,\reflectbox{\scriptsize\bf S}^{\frak X}} = \left\{ (x,X)\colon  x \in X, x\in\frak X , X\in\reflectbox{\bf S}^{\frak X}\right\} \subseteq \frak X  \times \reflectbox{\bf S}^{\frak X}
\end{equation}
as a \emph{membership relation} $x \in X$ between elements $x \in \frak X$ and subsets $X \in \reflectbox{S}^{\frak X}$; and also a complementary binary relation
\begin{equation}\label{0000-XsX-relation-complementary}
\mathscr{R}^c_{\frak X ,\reflectbox{\scriptsize\bf S}^{\frak X}} = \left\{ (x,X)\colon  x \not\in X, x\in\frak X , X\in\reflectbox{\bf S}^{\frak X}\right\} \subseteq \frak X  \times \reflectbox{\bf S}^{\frak X}
\end{equation}
as a \emph{non-membership relation} $x \not\in X$ between elements $x \in \frak X$ and subsets $X \in \reflectbox{S}^{\frak X}$; so that
\begin{equation}\label{0000-XsX-relations-sum}
\mathscr{R}_{\frak X ,\reflectbox{\scriptsize\bf S}^{\frak X}}+\mathscr{R}^c_{\frak X ,\reflectbox{\scriptsize\bf S}^{\frak X}}= \frak X  \times \reflectbox{\bf S}^{\frak X}.
\end{equation}
Finally, we add to the stock $\mbox{\boldmath$\Lambda$}$ so called \emph{terraced\footnote{Those who are familiar with the beginnings of the eventological theory \cite[2007]{Vorobyev2007} should keep their attention to the amazing inevitability of the ``splitting'' of the previously unified concept of the \emph{terrace event} into two dual halves, the right of which is the \emph{terraced ket-event} which is defined as a terrace event of the first kind $\displaystyle\textsf{ter}(X/\!\!/\frak X)=\bigcap_{x \in X} x \bigcap_{x \in\frak X -X} (\Omega-x)\subseteq\Omega$ from the eventological part of the Kolmogorov \emph{probability theory}, and the left one is a \emph{terraced bra-event}, a new concept from the \emph{theory of believabilities}, dual to the \emph{probability theory}, which is defined as terraced event of the 5th kind $\displaystyle\textsf{Ter}_{X/\!\!/\frak X}=\bigcup_{x \in X}
x\subseteq\Omega$ from the eventological classification.} label}
\begin{equation}\label{terraced}
\Big(\textsf{Ter}_{X/\!\!/\frak X}, \ \ \ \textsf{ter}(X/\!\!/\frak X)\Big) = \left( \bigcup_{x \in X} x, \ \ \ \bigcap_{x \in X} x \bigcap_{x \in\frak X -X} (\Omega-x) \right) \subseteq\Omega\times\Omega,
\end{equation}
numbered by labels-subsets $X\in\reflectbox{\bf S}^{\frak X}$ and while defined simply as a pair of indicated measurable subsets of
$\Omega$.

To have a full stock  we'll stock up in the literal sense ``complementary'' element-set labels, constructed from:
1) the complements
$x^c=\Omega-x$ to measurable subsets $x\subseteq\Omega$, 2) the \emph{Ì-complementary set} $\frak X ^{\!(\!c\!)}=\{x^c\colon
x\in\frak X \} \subseteq \mathcal{A}$ composed from these complements, and 3) the sets $\reflectbox{\bf
S}^{\frak X ^{\!(\!c\!)}}=\left\{X^{c(\!c\!)}\colon X\in \reflectbox{\bf
S}^{\frak X}\right\}\subseteq\mathscr{P}\left(\frak X ^{\!(\!c\!)}\right)$ of subsets
$X^{c(\!c\!)}=(X^c)^{\!(\!c\!)}=(\frak X -X)^{\!(\!c\!)}\subseteq\frak X ^{\!(\!c\!)}$, i.e., such that $X^{c(\!c\!)}=\{x^c\colon x\in X^c\} \in
\reflectbox{\bf S}^{\frak X ^{\!(\!c\!)}}$.

There we also place a label similar to (\ref{0000-XsX-cartesian}), the Cartesian product
\begin{equation}\label{0000-XsX-cartesian-c}
\frak X ^{\!(\!c\!)} \times \reflectbox{\bf S}^{\frak X ^{\!(\!c\!)}} = \left\{ (x^c,X^{c(\!c\!)}) : x^c \in \frak X ^{\!(\!c\!)}, X^{c(\!c\!)} \in \reflectbox{S}^{\frak X ^{\!(\!c\!)}} \right\},
\end{equation}
which defines analogous to (\ref{0000-XsX-relation-complementary}) a complementary binary relation
\begin{equation}\label{0000-XsX-relation-c}
\mathscr{R}^c_{\frak X ^{\!(\!c\!)},\reflectbox{\scriptsize\bf S}^{\frak X ^{\!(\!c\!)}}} = \left\{ (x^c,X^{c(\!c\!)})\colon  x^c \in X^{c(\!c\!)}, x^c\in\frak X ^{\!(\!c\!)}, X^{c(\!c\!)}\in\reflectbox{\bf S}^{\frak X ^{\!(\!c\!)}}\right\} \subseteq \frak X ^{\!(\!c\!)} \times \reflectbox{\bf S}^{\frak X ^{\!(\!c\!)}}
\end{equation}
as a \emph{membership relation} $x^c \in X^{c(\!c\!)}$ between elements $x^c \in \frak X ^{\!(\!c\!)}$ and subsets $X^{c(\!c\!)} \in
\reflectbox{S}^{\frak X ^{\!(\!c\!)}}$; and also a complementary binary relation
\begin{equation}\label{0000-XsX-relation-c-complementary-en}
\mathscr{R}_{\frak X ^{\!(\!c\!)},\reflectbox{\scriptsize\bf S}^{\frak X ^{\!(\!c\!)}}} = \left\{ (x^c,X^{c(\!c\!)})\colon  x^c \not\in X^{c(\!c\!)}, x^c\in\frak X ^{\!(\!c\!)}, X^{c(\!c\!)}\in\reflectbox{\bf S}^{\frak X ^{\!(\!c\!)}}\right\} \subseteq \frak X ^{\!(\!c\!)} \times \reflectbox{\bf S}^{\frak X ^{\!(\!c\!)}}
\end{equation}
as a \emph{non-membership relation} $x^c \not\in X^{c(\!c\!)}$ between elements $x^c \in \frak X ^{\!(\!c\!)}$ and subsets $X^{c(\!c\!)}
\in \reflectbox{S}^{\frak X ^{\!(\!c\!)}}$; so that
\begin{equation}\label{0000-XsX-relations-c-complementary-sum-en}
\mathscr{R}^c_{\frak X ^{\!(\!c\!)},\reflectbox{\scriptsize\bf S}^{\frak X ^{\!(\!c\!)}}}+
\mathscr{R}_{\frak X ^{\!(\!c\!)},\reflectbox{\scriptsize\bf S}^{\frak X ^{\!(\!c\!)}}}=
\frak X ^{\!(\!c\!)} \times \reflectbox{\bf S}^{\frak X ^{\!(\!c\!)}}.
\end{equation}
Finally, do not forget the similar to (\ref{terraced}) \emph{terrace label}
\begin{equation}\label{terraced-c}
\Big(\textsf{Ter}_{X^{\!c(\!c\!)}/\!\!/\frak X ^{\!(\!c\!)}}, \ \ \ \textsf{ter}\left(X^{c(\!c\!)}/\!\!/\frak X ^{\!(\!c\!)}\right)\Big) = \left( \bigcup_{x^c \in X^{\!c(\!c\!)}} x^c, \ \ \ \bigcap_{x^c \in X^{\!c(\!c\!)}} x^c \bigcap_{x^c \in\frak X ^{\!(\!c\!)}-X^{\!c(\!c\!)}} (\Omega-x^c) \right) \subseteq\Omega\times\Omega,
\end{equation}
numbered by labels-subsets $X^{\!c(\!c\!)}\in\reflectbox{\bf S}^{\frak X ^{\!(\!c\!)}}$.

The stock {\boldmath$\Lambda$} of element-set labels $\lambda\in\mbox{\boldmath$\Lambda$}$ is intended to construct such a system of element-set ``coordinates'', which, relying on a duality ``element--set'', will allow us to divide each concept of the \emph{theory of experience and chance (TEC)} into two dual parts and present it in the form of a conveniently written \emph{dual pair}, i.e., pairs composed of two dual parts. In the bra-ket notation \cite{Vorobyev2016famems1}, the dual parts of pairs labelled with the labels $\lambda,\lambda'\in\mbox{\boldmath$\Lambda$}$, are denoted by $\bra{\lambda}$ and $\ket{\lambda'}$ correspondingly, the entire dual pair is denoted by $\braket{\lambda|\lambda'}$ and is defined as the Cartesian product $\braket{\lambda|\lambda'}=\bra{\lambda}\times\ket{\lambda'}$ of their dual parts, placing the corresponding concept of the \emph{theory of experience and chance} in the system of ``element-set coordinates''.

\section{Co$\sim$event as a binary relation}%\\  \tiny$\sim$ \textsf{0000/vorobyev2-text, 28 April 2017}}

Let $\braket{\Omega,\mathcal A|\Omega,\mathcal A}=\big(\braket{\Omega|\Omega},\braket{\mathcal A|\mathcal A} \big)$ be a measurable \emph{bra-ket space}\footnote{In the following discourse, I use the notions and notations of the element-set labelling, introduced in the preliminary work \cite{Vorobyev2016famems1}, without necessarily defining them again here for the sake of space saving.}, labelled by the \emph{measurable binary relation} $\mathscr R\subseteq\braket{\Omega|\Omega}$ using $\mathscr R$-labels from the measurable space $(\Omega,\mathcal A)$ with $\mathscr R$-labelling sets $\frak X_{\!\mathscr R}\subseteq\mathcal A$ and $\reflectbox{\bf S}^{\frak X_{\!\mathscr R}}\subseteq\mathscr P(\frak X_{\!\mathscr R})$, which are defined the following way \cite{Vorobyev2016famems1}.

\texttt{\linebreak\indent Definition \!\refstepcounter{ctrdef}\arabic{ctrdef}\,\label{def-R-labelling-set}\itshape\footnotesize (basic  $\mathscr{R}$-labelling set $\frak X _\mathscr{R}$)\!.}
The \emph{basic $\mathscr{R}$-labelling set $\frak X _\mathscr{R}\subseteq\mathcal{A}$ of measurable subset of $\Omega$} is defined by the binary relation $\mathscr{R}\subseteq\braket{\Omega|\Omega}$ as the set of labels
\begin{equation}\label{R-labelling-X_R}
\frak X_{\!\mathscr R} = \left\{ x\in\mathcal A\colon \ket{x}=\mathscr R|_{\bra{\omega^*}}, \bra{\omega^*}\in\bra{\Omega}  \right\}\subseteq\mathcal{A},
\end{equation}
composed from measurable subsets $x\subseteq\Omega$ labelling ket-subsets $\ket{x}\subseteq\ket{\Omega}$ that serve by values of the cross-sections: $\ket{x}=\mathscr{R}|_{\bra{\omega^*}}\subseteq\ket{\Omega}$ of binary relation $\mathscr{R}$ by bra-points $\bra{\omega^*}\in\bra{\Omega}$.

Note, if there is the bra-point $\bra{\omega^*}\in\bra{\Omega}$ such that $\mathscr R|_{\bra{\omega^*}}=\varnothing_{\ket{\Omega}}$ then $\mathscr R|_{\bra{\omega^*}}=\ket{\varnothing^{\Omega}}$, i.e. the empty cross-section $\mathscr R|_{\bra{\omega^*}}$ coincides with the ket-subset $\ket{\varnothing^{\Omega}}$ where $\varnothing^{\Omega}\in\frak X_{\mathscr R}$ is the $\Omega$-empty label.

\texttt{\linebreak\indent Definition \!\refstepcounter{ctrdef}\arabic{ctrdef}\,\label{def-R-labelling-set-of-subsets}\itshape\footnotesize (basic set $\reflectbox{\bf S}^{\frak X_{\!\mathscr R}}$ of $\mathscr{R}$-labelling subsets)\!.}
The \emph{basic set $\reflectbox{\bf S}^{\frak X_{\!\mathscr R}}\subseteq\mathscr{P}(\frak X _\mathscr{R}-\{\varnothing^\Omega\})$ of $\mathscr{R}$-labelling subsets of measurable subsets of $\Omega$} is defined by the binary relation $\mathscr{R}\subseteq\braket{\Omega|\Omega}$ as the set of set-labels
\begin{equation}\label{R-labelling-S^X}
\reflectbox{\bf S}^{\frak X_{\!\mathscr R}} = \left\{ X\subseteq\frak X _\mathscr{R}-\{\varnothing^\Omega\}\colon \textsf{ter}(X/\!\!/\frak X _\mathscr{R})\ne\varnothing^\Omega \right\}\subseteq\mathscr{P}(\frak X _\mathscr{R}-\{\varnothing^\Omega\}),
\end{equation}
composed only from labelling subsets $X\subseteq\frak X _\mathscr{R}$ that do not contain the $\Omega$-empty label: $\varnothing^\Omega\notin X$, and number the $\Omega$-nonempty terraced labels: $\textsf{ter}(X/\!\!/\frak X _\mathscr{R})\ne\varnothing^\Omega$.

The measurable relation $\mathscr R$ generates the following element-set $\mathscr R$-labelling quotient-sets.
\begin{equation}\label{def-labelling3}
\bra{\Omega}\!/\mathscr{R}=\Bra{\frak X _\mathscr{R}}=\{\bra{x}\colon x\in\frak X _\mathscr{R} \}
\end{equation}
is the $\mathscr R$-labelling \emph{bra-quotient-set} $\bra{\Omega}\!/\mathscr{R}$ by the binary relation $\mathscr{R}\subseteq\braket{\Omega|\Omega}$,
under which the labels $x \in \frak X _\mathscr{R}$ of labelling set $\frak X _\mathscr{R}$ label all \emph{bra-subsets} $\bra{x}\in\bra{\Omega}\!/\mathscr{R}$ of the quotient-set $\bra{\Omega}\!/\mathscr{R}$;
\begin{equation}\label{def-labelling4}
\begin{split}
\ket{\Omega}\!\!/\mathscr{R}=\Ket{\reflectbox{\bf S}^{\frak X_{\!\mathscr R}}} &=\left\{\ket{\textsf{ter}(X/\!\!/\frak X _\mathscr{R})}\colon X\in\reflectbox{\bf S}^{\frak X_{\!\mathscr R}}\right\}
\end{split}
\end{equation}
is the  $\mathscr R$-labelling \emph{ket-quotient-set} $\ket{\Omega}\!\!/\mathscr{R}$ by the binary relation $\mathscr{R}\subseteq\braket{\Omega|\Omega}$,
under which the subsets $X \in \reflectbox{\bf S}^{\frak X_{\!\mathscr R}}$ from the set of labelling subsets $\reflectbox{\bf S}^{\frak X_{\!\mathscr R}}$ label the \emph{terraced ket-subsets} $\ket{\textsf{ter}(X/\!\!/\frak X _\mathscr{R})}\in\ket{\Omega}\!\!/\mathscr{R}$ of the quotient-set $\ket{\Omega}\!\!/\mathscr{R}$;
\begin{equation}\label{def-labelling5}
\begin{split}
\braket{\Omega|\Omega}\!\!/\mathscr{R}&=\Braket{\frak X _\mathscr{R}|\reflectbox{\bf S}^{\frak X_{\!\mathscr R}}}=\left\{\braket{x|\textsf{ter}(X/\!\!/\frak X _\mathscr{R})}\colon x\in\frak X _\mathscr{R},X\in\reflectbox{\bf S}^{\frak X_{\!\mathscr R}}\right\}
\end{split}
\end{equation}
is the $\mathscr R$-labelling \emph{bra-ket-quotient-set} $\braket{\Omega|\Omega}\!\!/\mathscr{R}$ by the binary relation $\mathscr{R}\subseteq\braket{\Omega|\Omega}$,
under which the pairs $(x,X)$, where $x \in \frak X _\mathscr{R}$ is an element of the labelling set $\frak X _\mathscr{R}$,
and $X \in \reflectbox{\bf S}^{\frak X_{\!\mathscr R}}$ is a subset from the set $\reflectbox{\bf S}^{\frak X_{\!\mathscr R}}$ of labelling subsets, label all \emph{bra-ket-subsets} $\braket{x|\textsf{ter}(X/\!\!/\frak X _\mathscr{R})}\in\braket{\Omega|\Omega}\!/\mathscr{R}$ of the quotient-set $\braket{\Omega|\Omega}\!\!/\mathscr{R}$.

\texttt{\linebreak\indent Predefinition \!\refstepcounter{ctrdefpre}\arabic{ctrdefpre}\,\label{defpre-all-events}\itshape\footnotesize (events and co$\sim$events)\!.}

$\star$ The bra-points $\bra{\omega}\in\ket{\Omega}$ are called  \emph{elementary bra-incomes (incomes)}.

$\star$ The bra-subsets $\bra{x}\subseteq\bra{\Omega}$ and terraced bra-subsets $\bra{\textsf{Ter}_{X/\!\!/\frak X_{\!\mathscr R}}}\subseteq\bra{\Omega}$ of the bra-set $\bra{\Omega}$ are called \emph{bra-events} and terraced \emph{bra-events} correspondingly.

$\star$ The ket-points $\ket{\omega}\in\ket{\Omega}$ are called \emph{elementary ket-outcomes (outcomes)}.

$\star$ The ket-subsets $\ket{x}\subseteq\ket{\Omega}$ and terraced ket-subsets $\ket{\textsf{ter}(X/\!\!/\frak X_{\!\mathscr R})}\subseteq\ket{\Omega}$ of the ket-set $\ket{\Omega}$ are called \emph{ket-events} and terraced \emph{ket-events} correspondingly.

$\star$ The bra-ket-subsets $\braket{x|x}\subseteq\braket{\Omega|\Omega}$, $\braket{\textsf{Ter}_{X/\!\!/\frak X_{\!\mathscr R}}|\textsf{ter}(X/\!\!/\frak X_{\!\mathscr R})}\subseteq\braket{\Omega|\Omega}$ and $\braket{x|\textsf{ter}(X/\!\!/\frak X_{\!\mathscr R})}\subseteq\braket{\Omega|\Omega}$ are called \emph{elementary bra-ket-events}.

$\star$ The bra-ket-subset $\mathscr R\subseteq\braket{\Omega|\Omega}$, i.e., any measurable binary relation, generating the $\mathscr R$-labelling, is called a \emph{co$\sim$event (an experienced-random co$\sim$event)}.

\texttt{\linebreak\indent Predefinition \!\refstepcounter{ctrdefpre}\arabic{ctrdefpre}\,\label{defpre-all-events}\itshape\footnotesize (full-believable, certainty, non-experienced, and impossible events and full-believable-certainty and non-experienced-impossible co$\sim$events)\!.}

$\star$ The \emph{bra-events} $\bra{\Omega}$ and $\bra{\varnothing}$ are called \emph{full-believable and non-experienced} correspondingly.

$\star$ The \emph{ket-events} $\ket{\Omega}$ and  $\ket{\varnothing}$ are called \emph{certainty and impossible} correspondingly.

$\star$ The \emph{co$\sim$events} $\braket{\Omega|\Omega}$  and $\braket{\varnothing|\varnothing}$ are called \emph{full-believable-certainty and non-experienced-impossible} correspondingly.

$\star$ The \emph{co$\sim$events} $\braket{\Omega|x}$ and $\braket{x|\Omega}$ are called \emph{full-believable-random and experienced-certainty} correspondingly.

$\star$ The \emph{co$\sim$events} $\braket{\varnothing|x}$ and  $\braket{x|\varnothing}$ are called \emph{non-experienced-random and experienced-impossible} correspondingly.

$\star$ The \emph{co$\sim$events} $\braket{\Omega|\varnothing}$ and $\braket{\varnothing|\Omega}$ are called \emph{full-believable-impossible and non-experienced-certainty} correspondingly.

\texttt{\linebreak\indent Predefinition \!\refstepcounter{ctrdefpre}\arabic{ctrdefpre}\,\label{defpre-R-labelled-events}\itshape\footnotesize  ($\mathscr{R}$-labelled events)\!.}
For the sake of brevity, the following general notation of \emph{$\mathscr{R}$-labelled events}, and suited general denotations:
\begin{equation}\label{R-labelled-events}
\begin{split}
\bra{\lambda^*_\mathscr{R}}&=
\begin{cases}
\bra{x}, & x\in \frak X _\mathscr{R},\\
\bra{\textsf{Ter}_{X/\!\!/\frak X _\mathscr{R}}}, & X\in \reflectbox{\bf S}^{\frak X_{\!\mathscr R}},
\end{cases}\\
&\\
\ket{\lambda_\mathscr{R}}&=
\begin{cases}
\ket{x}, & x\in \frak X _\mathscr{R},\\
\ket{\textsf{ter}(X/\!\!/\frak X _\mathscr{R})}, & X\in \reflectbox{\bf S}^{\frak X_{\!\mathscr R}},
\end{cases}\\
&\\
\braket{\lambda^*_\mathscr{R}|\lambda_\mathscr{R}}&=
\begin{cases}
\braket{x|x}, & x\in \frak X _\mathscr{R},\\
\braket{\textsf{Ter}_{X/\!\!/\frak X _\mathscr{R}}|\textsf{ter}(X/\!\!/\frak X _\mathscr{R})}, & X\in \reflectbox{\bf S}^{\frak X_{\!\mathscr R}},
\end{cases}\\
&\\
\braket{\lambda^*_\mathscr{R}|\lambda'_\mathscr{R}}&=
\braket{x|\textsf{ter}(X/\!\!/\frak X _\mathscr{R})}, \ x\in \frak X _\mathscr{R}, \, X\in \reflectbox{\bf S}^{\frak X_{\!\mathscr R}}
\end{split}
\end{equation}
are introduced for
\emph{ket-events} $\ket{x}\subseteq\ket{\Omega}$, \emph{terraced ket-events} $\ket{\textsf{ter}(X/\!\!/\frak X _\mathscr{R})}\subseteq\ket{\Omega}$, \emph{bra-events} $\bra{x}\subseteq\bra{\Omega}$, \emph{terraced bra-events} $\ket{\textsf{Ter}_{X/\!\!/\frak X _\mathscr{R}}}\subseteq\ket{\Omega}$, and also for \emph{elementary bra-ket-events}: $\braket{x|x}\!\subseteq\!\braket{\Omega|\Omega}$, $\braket{\textsf{Ter}_{X/\!\!/\frak X _\mathscr{R}}|\textsf{ter}(X/\!\!/\frak X _\mathscr{R})}\subseteq\braket{\Omega|\Omega}$ and $\braket{x|\textsf{ter}(X/\!\!/\frak X _\mathscr{R})}\subseteq\braket{\Omega|\Omega}$;
which are defined in Predefinition \ref{defpre-all-events} and labelled by the co$\sim$event $\mathscr{R}\subseteq\braket{\Omega|\Omega}$.

\texttt{\linebreak\indent Predefinition \!\refstepcounter{ctrdefpre}\arabic{ctrdefpre}\,\label{defpre-R-events-duality}\itshape\footnotesize (bra-ket-duality of $\mathscr{R}$-labelled events)\!.}
They say that the $\mathscr{R}$-labelled bra-event $\bra{\lambda^*_\mathscr{R}}$ and the $\mathscr{R}$-labelled ket-event $\ket{\lambda_\mathscr{R}}$ are \emph{bra-ket-dual each other} and form the \emph{pair of bra-ket-dual events} as the Cartesian product $\braket{\lambda^*_\mathscr{R}|\lambda_\mathscr{R}}=\bra{\lambda^*_\mathscr{R}}\times\ket{\lambda_\mathscr{R}}$.

%\clearpage
\section{``Something happens when that is experienced, what happens''}%\\  \tiny$\sim$ \textsf{0000/vorobyev2-text, May 2017}}

{\footnotesize
\hfill \emph{Die Welt ist alles, was der Fall ist.}\footnote{``The world is all that is the case.''}\\[-17pt]

}

{\scriptsize
\hfill Ludwig Wittgenstein \cite[1921]{Wittgenstein1921}
}\\[-16pt]

\smallskip

{\footnotesize
\hfill \emph{The world is all that occurs,}\\[-17pt]

\hfill \emph{when that is experienced, what happens.}\\[-17pt]

}

{\scriptsize
\hfill Theory of experience and chance [2017]
}\\[-16pt]

\subsection{The axiom of co$\sim$event as of what occurs, when that is experienced, what happens}%\\  \tiny$\sim$ \textsf{0000/vorobyev2-text, May 2017}}

Before the axioms \ref{Ax-co-event} (the axiom of co$\sim$event), which is central to the theory of experience and chance, I will formulate for comparison, in the same notation, what I called the ``silent'' Kolmogorov axiom. Its number is zero.

\texttt{\linebreak\indent Axiom \!\refstepcounter{ctrAx}\arabic{ctrAx}\,\label{Ax-event-Kolmogorov}\itshape\footnotesize (an event happens, when its elementary outcome happens [Kolmogorov theory of probabilities])\!.}

\texttt{(1)} \emph{The elementary outcome $\omega\in\Omega$ is what happens: $\omega=\omega^{\!\uparrow}$, or does not happen: $\omega\ne\omega^{\!\uparrow}$}.

\texttt{(2)} \emph{Any event  $\lambda\subseteq\Omega$ happens: $\lambda=\lambda^\uparrow$, when the elementary outcome happens: $\omega=\omega^\uparrow$, \textbf{which belong to it}: $\omega^\uparrow\in \lambda$}.

\texttt{\linebreak\indent Axiom \!\refstepcounter{ctrAx}\arabic{ctrAx}\,\label{Ax-co-event}\itshape\footnotesize (co$\sim$event occurs, when that is experienced, what happens [theory of experience and chance])\!.}

\texttt{(1)} \emph{The elementary ket-outcome $\ket{\omega}\in\ket{\Omega}$ is what happens: $\ket{\omega}=\ket{\omega}^{\!\uparrow}$, or does not happen: $\ket{\omega}\ne\ket{\omega}^{\!\uparrow}$}.

\texttt{(2)} \emph{For any $\mathscr R\subseteq\braket{\Omega|\Omega}$ any $\mathscr{R}$-labelled ket-event  $\ket{\lambda_{\!\mathscr R}}\subseteq\ket{\Omega}$ happens: $\ket{\lambda_{\!\mathscr R}}=\ket{\lambda_{\!\mathscr R}}^\uparrow$, when the elementary outcome happens: $\ket{\omega}=\ket{\omega}^\uparrow$, \textbf{which belong to it}: $\ket{\omega}^\uparrow\in \ket{\lambda_{\!\mathscr R}}$}.

\texttt{(3)} \emph{For any $\mathscr R\subseteq\braket{\Omega|\Omega}$ any $\mathscr{R}$-labelled bra-event  $\bra{\lambda^*_{\!\mathscr R}}\subseteq\bra{\Omega}$ is experienced: $\bra{\lambda^*_{\!\mathscr R}}=\bra{\lambda^*_{\!\mathscr R}}^\downarrow$, when dual $\mathscr{R}$-labelled ket-event happens: $\ket{\lambda_{\!\mathscr R}}=\ket{\lambda_{\!\mathscr R}}^\uparrow$}.

\texttt{(4)} \emph{The elementary bra-income $\bra{\omega^*}\in\bra{\Omega}$ is expereinced $\bra{\omega^*}=\bra{\omega^*}^\downarrow$, when  $\mathscr{R}$-labelled bra-event: $\bra{\lambda^*_{\!\mathscr R}}=\bra{\lambda^*_{\!\mathscr R}}^\downarrow$ is experienced, \textbf{to which $\bra{\omega^*}$ belongs}}: $\bra{\omega^*}\in\bra{\lambda^*_{\!\mathscr R}}^\downarrow$.

\texttt{(5)} \emph{The elementary income-outcome $\braket{\omega^*|\omega}\in\braket{\Omega|\Omega}$ is what occurs: $\braket{\omega^*|\omega}=\braket{\omega^*|\omega}^{\downarrow\uparrow}$, when the elementary ket-outcome: $\ket{\omega}=\ket{\omega}^{\!\uparrow}$ happens and the elementary bra-income: $\bra{\omega^*}=\bra{\omega^*}^\downarrow$ is experiencred; or does not occur: $\braket{\omega^*|\omega}\ne\braket{\omega^*|\omega}^{\downarrow\uparrow}$, when $\ket{\omega}\ne\ket{\omega}^{\!\uparrow}$ or $\bra{\omega^*}\ne\bra{\omega^*}^\downarrow$}.

\texttt{(6)} \emph{The co$\sim$event $\mathscr{R}\subseteq\braket{\Omega|\Omega}$ occurs: $\mathscr{R}=\mathscr{R}^{\downarrow\uparrow}$, when the elementary income-outcome: $\braket{\omega^*|\omega}=\braket{\omega^*|\omega}^{\downarrow\uparrow}$ occurs, \textbf{which belongs to it}}: $\braket{\omega^*|\omega}^{\downarrow\uparrow}\in\mathscr{R}$.

\subsection{Kolmogorov axioms}%\\  \tiny$\sim$ \textsf{0000/0001-axiom-ex-ch-Kolmogorov, 2 May 2017}}

\subsubsection{Kolmogorov axioms of believability theory}%\\  \tiny$\sim$ \textsf{0000/0001-axiom-ex-ch-Kolmogorov, 2 May 2017}}

Let $\bra{\Omega}$ be the bra-set of bra-points $\bra{\omega}\in\bra{\Omega}$, which we shall call the \emph{elementary bra-incomes} (or simply the \emph{elementary incomes}), and $\bra{\mathcal A}$ be the set of subsets from $\bra{\Omega}$. For any $\mathscr R\subseteq\braket{\Omega|\Omega}$ elements $\bra{\lambda^*_{\!\mathscr R}}\in\bra{\mathcal A}$ are called the $\mathscr R$-labelled \emph{bra-events}, and $\bra{\Omega}$ be the \emph{bra-set of elementary incomes}.

\texttt{\linebreak\indent Axiom \!\refstepcounter{ctrAx}\arabic{ctrAx}\,\label{Ax-bra-algebra}\itshape\footnotesize (algebra of bra-events)\!.} \emph{$\bra{\mathcal A}$ is an algebra of bra-events. The algebra of bra-events is also called the bra-algebra.}\\
\texttt{\linebreak\indent Axiom \!\refstepcounter{ctrAx}\arabic{ctrAx}\,\label{Ax-believability}\itshape\footnotesize (believability of bra-events)\!.} \emph{For any $\mathscr R\subseteq\braket{\Omega|\Omega}$ each $\mathscr R$-labelled bra-event $\bra{\lambda^*_{\!\mathscr R}}\in\bra{\mathcal A}$ is assigned the nonnegative real number $\mathbf B(\bra{\lambda^*_{\!\mathscr R}})$. This number is called the believability of $\mathscr R$-labelled bra-event $\bra{\lambda^*_{\!\mathscr R}}$.}\\
\texttt{\linebreak\indent Axiom \!\refstepcounter{ctrAx}\arabic{ctrAx}\,\label{Ax-believability-normalized}\itshape\footnotesize (normalization of believability)\!.} $\mathbf B(\bra{\Omega})=1$.\\
\texttt{\linebreak\indent Axiom \!\refstepcounter{ctrAx}\arabic{ctrAx}\,\label{Ax-believability-additivity}\itshape\footnotesize (additivity of believability)\!.} \emph{If $\mathscr R$-labelled bra-events $\bra{\lambda^*_{\!\mathscr R}}$ and $\bra{{\lambda^*_{\!\mathscr R}}'}$ are not intersected in $\bra{\Omega}$, then
\begin{equation}\label{believability-additiveness}\nonumber
\mathbf B(\bra{\lambda^*_{\!\mathscr R}}+\bra{{\lambda^*_{\!\mathscr R}}'})=\mathbf B(\bra{\lambda^*_{\!\mathscr R}})+\mathbf B(\bra{{\lambda^*_{\!\mathscr R}}'}).
\end{equation}}\\
\texttt{\linebreak\indent Axiom \!\refstepcounter{ctrAx}\arabic{ctrAx}\,\label{Ax-believability-continuity}\itshape\footnotesize (continuity of believability)\!.} \emph{For a decreasing sequence  $\bra{\lambda^*_{\!\mathscr R}}_1\supseteq\bra{\lambda^*_{\!\mathscr R}}_2\supseteq\ldots\supseteq\bra{\lambda^*_{\!\mathscr R}}_n\supseteq\ldots$ of $\mathscr R$-labelled bra-events from $\bra{\mathcal A}$ such that $\displaystyle\bigcap_n \bra{\lambda^*_{\!\mathscr R}}_n=\varnothing^{\bra{\Omega}}$, the equality $\displaystyle \lim_n \mathbf B(\bra{\lambda^*_{\!\mathscr R}}_n)=0$ takes place.}

Aggregate of objects $\bra{\Omega,\mathcal A,\mathbf B}=(\bra{\Omega},\bra{\mathcal A},\mathbf B)$, which is satisfied to axioms \ref{Ax-bra-algebra}, \ref{Ax-believability}, \ref{Ax-believability-normalized}, \ref{Ax-believability-additivity} and \ref{Ax-believability-continuity} we shall call the \emph{believability bra-space}, or simply the \emph{believability  space}.

\subsubsection{Kolmogorov axioms of probability theory}%\\  \tiny$\sim$ \textsf{0000/0001-axiom-ex-ch-Kolmogorov, 2 May 2017}}

Let $\ket{\Omega}$ be the ket-set of ket-points $\ket{\omega}\in\ket{\Omega}$ which we shall call the \emph{elementary ket-outcomes} (or simply the \emph{elementary outcomes}), and $\ket{\mathcal A}$ be the set of subsets from $\ket{\Omega}$. For any $\mathscr R\subseteq\braket{\Omega|\Omega}$ the elements $\ket{\lambda_{\!\mathscr R}}\in\ket{\mathcal A}$ of the set $\ket{\mathcal A}$ we shall call the \emph{$\mathscr R$-labelled ket-events}, and $\ket{\Omega}$ be the \emph{ket-set of elementary outcomes}.

\texttt{\linebreak\indent Axiom \!\refstepcounter{ctrAx}\arabic{ctrAx}\,\label{Ax-ket-algebra}\itshape\footnotesize (algebra of ket-events)\!.} \emph{$\ket{\mathcal A}$ is an algebra of ket-events. The algebra of ket-events is also called the ket-algebra.}\\
\texttt{\linebreak\indent Axiom \!\refstepcounter{ctrAx}\arabic{ctrAx}\,\label{Ax-probability}\itshape\footnotesize (probability of ket-events)\!.} \emph{For any $\mathscr R\subseteq\braket{\Omega|\Omega}$ each $\mathscr R$-labelled ket-event $\ket{\lambda_{\!\mathscr R}}\in\ket{\mathcal A}$ is assigned the nonnegative real number $\mathbf P(\ket{\lambda_{\!\mathscr R}})$. This number is called the probability of $\mathscr R$-labelled ket-event $\ket{\lambda_{\!\mathscr R}}$.}\\
\texttt{\linebreak\indent Axiom \!\refstepcounter{ctrAx}\arabic{ctrAx}\,\label{Ax-probability-normalized}\itshape\footnotesize (normalization of probability)\!.} $\mathbf P(\ket{\Omega})=1$.\\
\texttt{\linebreak\indent Axiom \!\refstepcounter{ctrAx}\arabic{ctrAx}\,\label{Ax-probability-additivity}\itshape\footnotesize (additivity of probability)\!.} \emph{If $\mathscr R$-labelled ket-events $\ket{\lambda_{\!\mathscr R}}$ and $\ket{\lambda_{\!\mathscr R}'}$ are not intersected in $\ket{\Omega}$, then}
\begin{equation}\label{probability-additiveness}\nonumber
\mathbf P(\ket{\lambda_{\!\mathscr R}}+\ket{\lambda_{\!\mathscr R}'})=\mathbf P(\ket{\lambda_{\!\mathscr R}})+\mathbf P(\ket{\lambda_{\!\mathscr R}'}).
\end{equation}

\texttt{\linebreak\indent Axiom \!\refstepcounter{ctrAx}\arabic{ctrAx}\,\label{Ax-probability-continuity}\itshape\footnotesize (continuity of probability)\!.} \emph{For a decreasing sequence  $\ket{\lambda_{\!\mathscr R}}_1\supseteq\ket{\lambda_{\!\mathscr R}}_2\supseteq\ldots\supseteq\ket{\lambda_{\!\mathscr R}}_n\supseteq\ldots$ of $\mathscr R$-labelled ket-events from $\ket{\mathcal A}$ such that $\displaystyle\bigcap_n \ket{\lambda_{\!\mathscr R}}_n=\varnothing^{\ket{\Omega}}$ the equality $\displaystyle \lim_n \mathbf P(\ket{\lambda_{\!\mathscr R}}_n)=0$ takes place.}

The aggregate of oblects $\ket{\Omega,\mathcal A,\mathbf P}=(\ket{\Omega},\ket{\mathcal A},\mathbf P)$, which is satisfied to axioms \ref{Ax-ket-algebra}, \ref{Ax-probability}, \ref{Ax-probability-normalized}, \ref{Ax-probability-additivity} and \ref{Ax-probability-continuity} we shall call the \emph{probability ket-space}, or simply the \emph{probability space}.

\subsection{Axioms of the theory of certainties (believabilities-probabilities)}%\\  \tiny$\sim$ \textsf{0000/0001-axiom-ex-ch-Kolmogorov, 3 May 2017}}

Let $\braket{\Omega|\Omega}=\bra{\Omega}\times\ket{\Omega}$ be the set of bra-ket-points $\braket{\omega^*|\omega}=\bra{\omega^*}\times\ket{\omega}\in\braket{\Omega|\Omega}$, which we shall call  the \emph{elementary bra-ket-incomes-outcomes} (or simply the \emph{elementary incomes-outcomes}), and $\braket{\mathcal A|\mathcal A}$ be the set of subsets from $\braket{\Omega|\Omega}$. For any $\mathscr R\subseteq\braket{\Omega|\Omega}$ the elements $\braket{\lambda^*_{\!\mathscr R}|\lambda_{\!\mathscr R}}\in\braket{\mathcal A|\mathcal A}$ are called the \emph{$\mathscr R$-labelled bra-ket-events}, and $\braket{\Omega|\Omega}$ be the \emph{bra-ket-set of elementary incomes-outcomes}.

\texttt{\linebreak\indent Axiom \!\refstepcounter{ctrAx}\arabic{ctrAx}\,\label{Ax-bra-ket-algebra}\itshape\footnotesize (algebra of bra-ket-events)\!.} \emph{$\braket{\mathcal A|\mathcal A}=\alpha\big(\bra{\mathscr A}\times\ket{\mathscr A}\big)$ is a minimal algebra of bra-ket-events, which contains the Cartesian product of algebras $\bra{\mathscr A}\times\ket{\mathscr A}$. This algebra is also called bra-ket-algebra.}\\
\texttt{\linebreak\indent Axiom \!\refstepcounter{ctrAx}\arabic{ctrAx}\,\label{Ax-certainty}\itshape\footnotesize (certainty of bra-ket-events)\!.} \emph{For any $\mathscr R\subseteq\braket{\Omega|\Omega}$ each $\mathscr R$-labelled bra-ket-event $\braket{\lambda^*_{\!\mathscr R}|\lambda_{\!\mathscr R}}\in\braket{\mathcal A|\mathcal A}$ is assigned the nonnegative real number $\bfPhi(\braket{\lambda^*_{\!\mathscr R}|\lambda_{\!\mathscr R}})=\mathbf B(\bra{\lambda_{\!\mathscr R}^*})\mathbf P(\ket{\lambda_{\!\mathscr R}})$. This number is called the certainty of $\mathscr R$-labelled bra-ket-event $\braket{\lambda^*_{\!\mathscr R}|\lambda_{\!\mathscr R}}$.}\\
\texttt{\linebreak\indent Property \!\refstepcounter{ctrATT}\arabic{ctrATT}\,\label{ATT-certainty-normalized}\itshape\footnotesize (normalization of certainty)\!.} $\bfPhi(\braket{\Omega|\Omega})=1$.\\
\texttt{Proof.} $\bfPhi(\braket{\Omega|\Omega})=\mathbf B(\bra{\Omega})\mathbf P(\ket{\Omega})=1$ by axioms \ref{Ax-believability-normalized}, \ref{Ax-probability-normalized} and \ref{Ax-certainty}.\\
\texttt{\linebreak\indent Property \!\refstepcounter{ctrATT}\arabic{ctrATT}\,\label{ATT-certainty-additivity}\itshape\footnotesize (additivity of certainty)\!.} \emph{If $\mathscr R$-labelled bra-ket-events $\braket{\lambda^*_{\!\mathscr R}|\lambda_{\!\mathscr R}}$ and $\braket{{{\lambda_{\!\mathscr R}^*}}'|\lambda_{\!\mathscr R}'}$ are not intersected in $\braket{\Omega|\Omega}$, then}
\begin{equation}\label{certainty-additivity}\nonumber
\bfPhi(\braket{\lambda^*_{\!\mathscr R}|\lambda_{\!\mathscr R}}+\braket{{{\lambda_{\!\mathscr R}^*}}'|\lambda_{\!\mathscr R}'})=\bfPhi(\braket{\lambda^*_{\!\mathscr R}|\lambda_{\!\mathscr R}})+\bfPhi(\braket{{{\lambda_{\!\mathscr R}^*}}'|\lambda_{\!\mathscr R}'}).
\end{equation}
\texttt{Proof.} An additivity of product of additive functions $\mathbf B$ and $\mathbf P$ is a fact of general measure theory. So that an additivity of certainty $\bfPhi$ on $\braket{\mathscr A|\mathscr A}$ follows routinely from the axioms \ref{Ax-believability-additivity}, \ref{Ax-probability-additivity} and \ref{Ax-certainty}.

\texttt{\linebreak\indent Property \!\refstepcounter{ctrATT}\arabic{ctrATT}\,\label{ATT-certainty-continuity}\itshape\footnotesize (continuity of certainty)\!.} \emph{For a decreasing sequence $\braket{\lambda^*_{\!\mathscr R}|\lambda_{\!\mathscr R}}_1\supseteq\braket{\lambda^*_{\!\mathscr R}|\lambda_{\!\mathscr R}}_2\supseteq\ldots\supseteq\braket{\lambda^*_{\!\mathscr R}|\lambda_{\!\mathscr R}}_n\supseteq\ldots$ of $\mathscr R$-labelled bra-ket-events from $\braket{\mathcal A|\mathcal A}$ such that $\displaystyle\bigcap_n \braket{\lambda^*_{\!\mathscr R}|\lambda_{\!\mathscr R}}_n=\varnothing^{\braket{\Omega|\Omega}}$, the equality $\displaystyle \lim_n \bfPhi(\braket{\lambda^*_{\!\mathscr R}|\lambda_{\!\mathscr R}}_n)=0$ takes place.}\\
\texttt{Proof.} A continuity of certainty $\bfPhi$ on $\braket{\mathcal A|\mathcal A}$ as a product of the continuous believability $\mathbf B$ on  $\bra{\mathcal A}$ and the continuous probability $\mathbf P$ on $\ket{\mathcal A}$ follows from the general measure theory by axioms \ref{Ax-believability-continuity}, \ref{Ax-probability-continuity} and \ref{Ax-certainty}.\\

The aggregate of objects $\braket{\Omega,\mathcal A,\mathbf B|\Omega,\mathcal A,\mathbf P}=(\braket{\Omega|\Omega},\braket{\mathcal A|\mathcal A},\bfPhi)$, which is satisfied to axioms \ref{Ax-bra-ket-algebra} and \ref{Ax-certainty} we shall call the \emph{certainty (believability-probability) bra-ket-space}, or simply the \emph{certainty space}.

\texttt{\linebreak\indent Property \!\refstepcounter{ctrATT}\arabic{ctrATT}\,\label{ATT-others-properties-}\itshape\footnotesize (believability, probability and certainty of some events and co$\sim$event)\!.} From the axioms of the theory of certainties it follows that

$\star$ $\bfPhi(\braket{\Omega|\Omega})=\mathbf B(\bra{\Omega})=\mathbf P(\ket{\Omega})=1$,

$\star$ $\bfPhi(\braket{\varnothing|\varnothing})=\mathbf B(\bra{\varnothing})=\mathbf P(\ket{\varnothing})=0$,

$\star$ $\bfPhi(\braket{\Omega|\varnothing})=\bfPhi(\braket{\varnothing|\Omega})=0$,

$\star$ $\bfPhi(\braket{\varnothing|\lambda_{\!\mathscr R}})=\bfPhi(\braket{\lambda^*_{\!\mathscr R}|\varnothing})=0$,

$\star$ $\bfPhi(\braket{\Omega|\lambda_{\!\mathscr R}})=\mathbf P(\ket{\lambda_{\!\mathscr R}})$,

$\star$ $\bfPhi(\braket{\lambda^*_{\!\mathscr R}|\Omega})=\mathbf B(\bra{\lambda^*_{\!\mathscr R}})$.

\texttt{\linebreak\indent Note \!\refstepcounter{ctrnot}\arabic{ctrnot}\,\label{not-infinity-spaces}\itshape\footnotesize (on infinity spaces)\!.} For an exhaustive presentation of innovations in the postulating of the theory of experience and the chance, it is quite sufficient to have the finite space and four first Kolmogorov axioms in both the believability bra-space (axioms \ref{Ax-bra-algebra} --- \ref{Ax-believability-additivity}), and the probability ket-space (axioms \ref{Ax-ket-algebra} --- \ref{Ax-probability-additivity}).
And to postulate the theory of experience and chance in infinite spaces, it takes only a long-known necessary, but routine procedure, to linger here on which I do not see any special need.
Therefore, dropping the routine, we will always assume that we have at our disposal the smallest sigma-algebras $\bra{\mathcal A}$, $\ket{\mathcal A}$ and $\braket{\mathcal A|\mathcal A}$, containing those sigma algebras that are sufficient for finite space; and the believability $\mathbf B$, the probability $\mathbf P$ and the certainty $\bfPhi$ are
countably additive functions obtained as a result of unique extensions to all sets from the corresponding sigma-algebras $\bra{\mathcal A}$, $\ket{\mathcal A}$ and $\braket{\mathcal A|\mathcal A}$. Thus, it is always assumed that the believability bra-space $\bra{\Omega,\mathcal A,\mathbf B}$, the probability ket-space $\ket{\Omega,\mathcal A,\mathbf P}$ and the certainty bra-ket-space $\braket{\Omega,\mathcal A,\mathbf B|\Omega,\mathcal A,\mathbf P}$ are \emph{Borel spaces}, so that the new theory of experience and chance had complete freedom of action, not connected with the danger of coming to events or to co$\sim$events, which have no believability, no probability or no certainty.

\subsection{Properties of co$\sim$events and its dual halves: bra-events and ket-events}%\\  \tiny$\sim$ \textsf{0000/0001-axiom-ex-ch-properties, 23 March 2017}}

\texttt{\linebreak\indent Property \!\refstepcounter{ctrATT}\arabic{ctrATT}\,\label{ATT-bra-ket-events1}\itshape\footnotesize (bra-event is experienced, when ket-event happens)\!.}
\emph{If the ket-event $\ket{x}\subseteq\ket{\Omega}$ happens: $\ket{x}=\ket{x}^{\!\uparrow}$, then the bra-event $\bra{x}\subseteq{\Omega}$ is experienced: $\bra{x}=\bra{x}^{\downarrow}$. Otherwise, when ket-event $\ket{x}\subseteq\ket{\Omega}$ does not happen: $\ket{x}\ne\ket{x}^{\!\uparrow}$, the bra-event isn't experienced: $\bra{x}\ne\bra{x}^{\downarrow}$}.

\texttt{Proof} follows from the item \texttt{(3)} of Axiom \ref{Ax-co-event}.

\texttt{\linebreak\indent Property \!\refstepcounter{ctrATT}\arabic{ctrATT}\,\label{ATT-bra-ket-events2}\itshape\footnotesize (bra-events from which something follows; ket-events that follow from something)\!.}

\texttt{(1)} \emph{If the ket-event} $\ket{x}\subseteq\ket{\Omega}$ \emph{happens:} $\ket{x}=\ket{x}^{\!\uparrow}$, \emph{then all ket-events which contain it:} $\ket{x}^{\!\uparrow}\subseteq\ket{y}\subseteq\ket{\Omega}$ \emph{happens}: $\ket{y}=\ket{y}^{\!\uparrow}$; \emph{in other words, all ket-events, which follow from $\ket{x}^{\!\uparrow}$, happens}.

\texttt{(2)} \emph{If the bra-event $\bra{x}\subseteq\bra{\Omega}$ is experienced: $\bra{x}=\bra{x}^\downarrow$, then all bra-events, which are contained in it: $\bra{y}\subseteq\bra{x}^\downarrow\subseteq\bra{\Omega}$, are experienced: $\bra{y}=\bra{y}^{\downarrow}$; in other words, all bra events, from which $\bra{x}^\downarrow$ follows, are experienced.}

\texttt{Proof} follows from the items (2) and (4) of Axiom \ref{Ax-co-event}.

\texttt{\linebreak\indent Property \!\refstepcounter{ctrATT}\arabic{ctrATT}\,\label{ATT-terraced-bra-ket-events}\itshape\footnotesize (terraced bra-event is experienced, terraced ket-event happens)\!.}

\texttt{(1)} \emph{The terraced ket-event}
$$
\ket{\textsf{ter}(X/\!\!/\frak X _\mathscr{R})}
=\bigcap_{x\in X}\ket{x} \bigcap_{x\in \frak X _\mathscr{R}-X}(\ket{\Omega}-\ket{x}) \in\ket{\mathscr{A}}
$$
\emph{happens:} $\ket{\textsf{ter}(X/\!\!/\frak X _\mathscr{R})}=\ket{\textsf{ter}(X/\!\!/\frak X _\mathscr{R})}^{\!\uparrow}$, \emph{when  the ket-outcome, which belongs to it:} $\ket{\omega}\in\ket{\textsf{ter}(X/\!\!/\frak X _\mathscr{R})}$, \emph{happens:}
$\ket{\omega}=\ket{\omega}^{\!\uparrow}$. \emph{Otherwise, the terraced ket-event does not happens:} $\ket{\textsf{ter}(X/\!\!/\frak X _\mathscr{R})}\ne\ket{\textsf{ter}(X/\!\!/\frak X _\mathscr{R})}^{\!\uparrow}$.

\texttt{(2)} \emph{The terraced bra-event}
$$
\displaystyle\bra{\textsf{Ter}_{X/\!\!/\frak X _\mathscr{R}}}=\sum_{x\in X}\bra{x}\in\bra{\mathscr{A}}
$$
\emph{is experienced:} $\bra{\textsf{Ter}_{X/\!\!/\frak X _\mathscr{R}}}=\bra{\textsf{Ter}_{X/\!\!/\frak X _\mathscr{R}}}^\downarrow$, \emph{when the terraced ket-event}
$\ket{\textsf{ter}(X'/\!\!/\frak X _\mathscr{R})}=\ket{\textsf{ter}(X'/\!\!/\frak X _\mathscr{R})}^{\!\uparrow}$, \emph{such that}
$X\subseteq X'$ (see Footnote\footnote{In the evenology \cite{Vorobyev2007} this event has a special denotation: $\displaystyle\textsf{ter}_{X/\!\!/\frak X_{\!\mathscr R}}=\sum_{X\subseteq X'}\textsf{ter}(X'/\!\!/\frak X_{\!\mathscr R})$ and is called the \emph{terraced event of the 2-d type.}}) \emph{happens. Otherwise, the terraced bra-event isn't experienced:} $\bra{\textsf{Ter}_{X/\!\!/\frak X _\mathscr{R}}}\ne\bra{\textsf{Ter}_{X/\!\!/\frak X _\mathscr{R}}}^\downarrow$.

\texttt{Proof} of \texttt{(1)} follows from the item \texttt{(2)} of Axiom \ref{Ax-co-event}, and the proof of \texttt{(2)}
follows from the item \texttt{(3)} of Axiom \ref{Ax-co-event} and the item \texttt{(1)} of Property \ref{ATT-bra-ket-events2}.

\texttt{\linebreak\indent Property \!\refstepcounter{ctrATT}\arabic{ctrATT}\,\label{ATT-bra-ket-events4}\itshape\footnotesize (co$\sim$event as a membership relation)\!.}
\emph{Any co$\sim$event $\mathscr{R}\subseteq\braket{\Omega|\Omega}$ in the measurable bra-ket-space $\braket{\Omega,\mathcal A|\Omega,\mathcal A}$ is equivalence to the membership relation}
\begin{equation}\label{R-membership}
\mathscr R_{\braket{\frak X _\mathscr{R}|\reflectbox{\scriptsize\bf S}^{\frak X_{\!\mathscr R}}}}=\Big\{
\braket{x|\textsf{ter}(X/\!\!/\frak X _\mathscr{R})}\colon x\in X
\Big\}
\subseteq\Braket{\frak X _\mathscr{R}|\reflectbox{\bf S}^{\frak X_{\!\mathscr R}}}
\end{equation}
\emph{on element-set $\mathscr{R}$-labelling $\braket{\frak X _\mathscr{R}|\reflectbox{\bf S}^{\frak X_{\!\mathscr R}}}$ of the quotient-set $\braket{\Omega|\Omega}\!\!/\mathscr{R}$. In other words,}
\begin{equation}\label{R-membership}
\mathscr R=\Big\{\braket{\omega^*|\omega}\in\braket{\Omega|\Omega}\colon \braket{\omega^*|\omega}\in\braket{x|\textsf{ter}(X/\!\!/\frak X _\mathscr{R})}\in\mathscr R_{\braket{\frak X _\mathscr{R}|\reflectbox{\scriptsize\bf S}^{\frak X_{\!\mathscr R}}}}
\Big\} \subseteq\braket{\Omega|\Omega}.
\end{equation}
\emph{Wherein the co$\sim$event $\mathscr R$ occurs then and only then, when the elementary income-outcome $\braket{\omega^*|\omega}=\braket{\omega^*|\omega}^{\downarrow\uparrow}$ occurs, such that $\braket{\omega^*|\omega}^{\downarrow\uparrow}\in\braket{x|\textsf{ter}(X/\!\!/\frak X _\mathscr{R})}$, and the membership relation: $x\in X$ holds.}

\texttt{Proof} relies on equivalence of the inclusion relation  $\subseteq_{\braket{\Omega|\Omega}}$ and the membership relation $\in_{\frak X_{\!\mathscr R}}$ (see \cite{Vorobyev2016famems1})
\begin{eqnarray}
\label{membership-relation2}
\braket{x|\textsf{ter}(X/\!\!/\frak X _\mathscr{R})}\subseteq\mathscr R \ &\Longleftrightarrow& \ x\in X,
\end{eqnarray}
from which it follows that the co$\sim$event $\mathscr R$ occurs, i.e., $\braket{\omega^*|\omega}^{\downarrow\uparrow}\in\mathscr R$, then and only then, when two membership relations
$\braket{\omega^*|\omega}^{\downarrow\uparrow}\in\braket{x|\textsf{ter}(X/\!\!/\frak X _\mathscr{R})}$ and $x\in X$ hold. This proves the property.

\section{Beliavability, probability and certainty (believability-probability) measures in the theory of experience and chance}%\\  \tiny$\sim$ \textsf{0000/0002-measures, 8 May 2017}}

For convenience, we introduce abbreviated notation for the probability, believability and certainty of some bra-events, ket-events and bra-ket-events\footnote{
The bra-events from the set $\bra{\frak X_{\!\mathscr R}}=\{\bra{x}\colon x\in\frak X_{\!\mathscr R}\}$ are disjoint and form a partition of the bra-space $\bra{\Omega}=\sum_{x\in\frak X_{\!\mathscr R}}\bra{x}$. The such set $\bra{\frak X_{\!\mathscr R}}$ it generates only two kinds (of the six standard kinds \cite{Vorobyev2007}) ``non-trivial'' terraced bra-events. These are terraced bra-events of the fifth kind $\bra{\textsf{Ter}_{\!X/\!\!/\frak X_{\!\mathscr R}}}=\sum_{x\in X}\bra{x}$ and of the third kind $\bra{\textsf{ter}^{X^{\!c}\!/\!\!/\frak X_{\!\mathscr R}}}=\bigcap_{x\in X^{\!c}}\bra{x}$, which, moreover, are corresponding complements of each other in the bra-space $\bra{\Omega}$:
$
\bra{\Omega}-\bra{\textsf{Ter}_{\!X/\!\!/\frak X_{\!\mathscr R}}} = \bra{\textsf{ter}_{\!X^{\!c}\!/\!\!/\frak X_{\!\mathscr R}}}, \ X\subseteq\frak X_{\!\mathscr R}
$.
The remaining four kinds of events are constants for all $X\subseteq\frak X_{\!\mathscr R}$, which are equal to either $\emptyset_{\bra{\Omega}}$, or$\bra{\Omega}$. This is easy to verify.
In the bra-ket formalism of the theory of experience and chance, one pair of dual terraces is singled out:
the terraced bra-event of the fifth kind $\bra{\textsf{Ter}_{\!X/\!\!/\frak X_{\!\mathscr R}}}$ and the terraced bra-event of the first kind $\ket{\textsf{ter}(X/\!\!/\frak X_{\!\mathscr R})}$, which serve as the dual ``vis-a-vis''.
And the terraced events of other kinds do not play a significant role. This fact and concern for the brevity of the bra-ket formalism is the reason for my deviation from the standard terraced designations: to indicate the believability of terraced bra-events of the fifth kind I use not the standard abbreviation $b_{X/\!\!/\frak X_{\!\mathscr R}}=\mathbf B(\bra{\textsf{Ter}_{\!X/\!\!/\frak X_{\!\mathscr R}}})$, the subindex in which repeats the subindex in the record of the terrace bra-event, but the abbreviation
$
b(X/\!\!/\frak X_{\!\mathscr R})=\mathbf B(\bra{\textsf{Ter}_{\!X/\!\!/\frak X_{\!\mathscr R}}})
$,
which clearly corresponds to the designation for the probability of its dual ``vis-a-vis'', the terraced ket-event of the first kind:
$
p(X/\!\!/\frak X_{\!\mathscr R})=\mathbf P(\bra{\textsf{ter}(X/\!\!/\frak X_{\!\mathscr R})})
$.}:
\begin{equation}\label{p_xb_xp(X)b(X)}
\begin{split}
b_x  = \mathbf B(\bra{x})
&\mbox{ \small --- believability of the bra-event } \bra{x}\in\bra{\mathcal A},\\
p_x = \mathbf P(\ket{x})
&\mbox{ \small --- probability of the ket-event } \ket{x}\in\ket{\mathcal A},\\
b(X/\!\!/\frak X_{\!\mathscr R}) = \mathbf B(\bra{\textsf{Ter}_{X/\!\!/\frak X_{\!\mathscr R}}})
&\mbox{ \small --- believability of the terraced bra-events } \bra{\textsf{Ter}_{X/\!\!/\frak X_{\!\mathscr R}}}\in\bra{\mathcal A},\\
p(X/\!\!/\frak X_{\!\mathscr R}) = \mathbf P(\ket{\textsf{ter}(X/\!\!/\frak X_{\!\mathscr R})})
&\mbox{ \small --- probability of the terraced ket-event } \ket{\textsf{ter}(X/\!\!/\frak X_{\!\mathscr R})}\in\ket{\mathcal A},\\
\bfphi_x=\bfPhi(\braket{x|x})
&\mbox{ \small --- certainty of the bra-ket-event } \braket{x|x}\in\braket{\mathcal A|\mathcal A},\\
\bfphi(X/\!\!/\frak X_{\!\mathscr R})=\bfPhi(\braket{\textsf{Ter}_{X/\!\!/\frak X_{\!\mathscr R}}|\textsf{ter}(X/\!\!/\frak X_{\!\mathscr R})})
&\mbox{ \small --- certainty of the bra-ket-event } \braket{\textsf{Ter}_{X/\!\!/\frak X_{\!\mathscr R}}|\textsf{ter}(X/\!\!/\frak X_{\!\mathscr R})}\in\braket{\mathcal A|\mathcal A},\\
\bfphi_x(X/\!\!/\frak X_{\!\mathscr R})=\bfPhi(\braket{x|\textsf{ter}(X/\!\!/\frak X_{\!\mathscr R})})
&\mbox{ \small --- certainty of the bra-ket-event } \braket{x|\textsf{ter}(X/\!\!/\frak X_{\!\mathscr R})}\in\braket{\mathcal A|\mathcal A}.
\end{split}\hspace{-30pt}
\end{equation}
By Axiom \ref{Ax-certainty} we have
\begin{equation}\label{p_xb_xp(X)b(X)}
\begin{split}
\bfphi_x=b_xp_x
&\mbox{ \small --- certainty of the bra-ket-event } \braket{x|x}\in\braket{\mathcal A|\mathcal A},\\
\bfphi(X/\!\!/\frak X_{\!\mathscr R})=b(X/\!\!/\frak X_{\!\mathscr R})p(X/\!\!/\frak X_{\!\mathscr R})
&\mbox{ \small --- certainty of the bra-ket-event } \braket{\textsf{Ter}_{X/\!\!/\frak X_{\!\mathscr R}}|\textsf{ter}(X/\!\!/\frak X_{\!\mathscr R})}\in\braket{\mathcal A|\mathcal A},\\
\bfphi_x(X/\!\!/\frak X_{\!\mathscr R})=b_xp(X/\!\!/\frak X_{\!\mathscr R})
&\mbox{ \small --- certainty of the bra-ket-event } \braket{x|\textsf{ter}(X/\!\!/\frak X_{\!\mathscr R})}\in\braket{\mathcal A|\mathcal A}.
\end{split}
\end{equation}

\texttt{\linebreak\indent Theorem \!\refstepcounter{ctrTh}\arabic{ctrTh}\,\label{Th-Fubini}\itshape\footnotesize (certainty of a co$\sim$event, Robbins-Fubini theorem \cite{Robbins1944/45, Fubini1907})\!.}
\emph{The certainty (believability-probability) $\bfPhi(\mathscr{R})=\bfPhi(\braket{\omega^*|\omega}\in\mathscr{R})$ of the co$\sim$event $\mathscr{R}\subseteq\braket{\Omega|\Omega}$ can be calculated from two equivalent formulas:}
\begin{eqnarray}
\label{fubini1}
\bfPhi(\mathscr{R})
&=&\sum_{x\in\frak X _\mathscr{R}} \bfphi_x,\\
\label{fubini2}
\bfPhi(\mathscr{R})
&=&\sum_{X\in\,\reflectbox{\scriptsize\bf S}^{\frak X_{\!\mathscr R}}} \bfphi(X/\!\!/\frak X _\mathscr{R}).
\end{eqnarray}

\texttt{Proof} of formulas (\ref{fubini1}) and (\ref{fubini2}) is based on a change in the order of the iterated sums and is analogous to the proof of the well-known theorem of Fubini on reducing the calculation of the double sum to the calculation of iterated sums:
\begin{equation}\label{fubini3}
\begin{split}
\bfPhi(\mathscr{R})
&=\sum_{x\in\frak X _\mathscr{R}}\sum_{x\in X\in\,\reflectbox{\scriptsize\bf S}^{\frak X_{\!\mathscr R}}}
\bfPhi\Big(\braket{\omega^*|\omega}\in\braket{x|\textsf{ter}(X/\!\!/\frak X _\mathscr{R})}\Big)\\
&=\sum_{x\in\frak X _\mathscr{R}}\sum_{x\in X\in\,\reflectbox{\scriptsize\bf S}^{\frak X_{\!\mathscr R}}} b_x p(X/\!\!/\frak X _\mathscr{R})
=\sum_{x\in\frak X _\mathscr{R}} b_x p_x
=\sum_{x\in\frak X _\mathscr{R}} \bfphi_x,
\end{split}
\end{equation}
\begin{equation}\label{fubini4}
\begin{split}
\bfPhi(\mathscr{R})
&=\sum_{X\in\,\reflectbox{\scriptsize\bf S}^{\frak X_{\!\mathscr R}}}
  \sum_{x\in X\in\,\reflectbox{\scriptsize\bf S}^{\frak X_{\!\mathscr R}}}
\bfPhi\Big(\braket{\omega^*|\omega}\in\braket{x|\textsf{ter}(X/\!\!/\frak X _\mathscr{R})}\Big)\\
&=\sum_{X\in\,\reflectbox{\scriptsize\bf S}^{\frak X_{\!\mathscr R}}}
  \sum_{x\in X\in\,\reflectbox{\scriptsize\bf S}^{\frak X_{\!\mathscr R}}}
  b_x p(X/\!\!/\frak X _\mathscr{R})
=\sum_{X\in\,\reflectbox{\scriptsize\bf S}^{\frak X_{\!\mathscr R}}}
  b(X/\!\!/\frak X _\mathscr{R}) p(X/\!\!/\frak X _\mathscr{R})
=\sum_{X\in\,\reflectbox{\scriptsize\bf S}^{\frak X_{\!\mathscr R}}}
  \bfphi(X/\!\!/\frak X _\mathscr{R}).
\end{split}
\end{equation}

\section{Experienced, random and experienced-random variables in the theory of experience and chance}%\\ \tiny $\sim$ \textsf{experienced-and-random-variables, 6 Apr 2016 $\sim$ !!!!!!}}

\emph{Experienced, random and experienced-random variables} are a part of the basic concepts of the \emph{theory of experience and chance}.
Complete and free from any unnecessary restrictions the presentation of the foundations of the theory of probabilities on the basis of measure theory is given by Kolmogorov \cite[1933]{Kolmogorov1933};
it made it quite obvious that the \emph{random variable} is nothing more than a measurable function on the \emph{probability space}.
The theory of experience and chance
also relies on the measure theory, which makes it equally obvious that the \emph{experienced variable} dual to random one, in turn, is nothing more than a measurable function on the \emph{believability space} dual to the probability one.
An experienced-random variable is defined as a measurable function on the Cartesian product of believability and probability spaces, the \emph{certainty space}.
These circumstances can not be ignored in the presentation of the beginning of the theory of experience and chance, which succeeded in combining the theory of believabilities and the theory of probabilities on the basis of the concepts of the \emph{space of elementary incomes} and the \emph{space of elementatry outcomes} and their Cartesian product, the \emph{space of elementary incomes-outsomes}, and one must not forget, each time emphasizing, that only when one is immersed in a dual context of the theory of experience and chance, representations about \emph{experienced, random and experienced-random variables} acquire the mathematical and applied content.

\subsection{Experienced variable}%\\ \tiny $\sim$ \textsf{experienced-and-random-variables, 10 May 2017}}

\texttt{\linebreak\indent Definition \!\refstepcounter{ctrdef}\arabic{ctrdef}\,\label{def-experienced-variable}\itshape\footnotesize (experienced variable)\!.} The function ${\bra{\xi_{\!\mathscr R}}}:\bra{\Omega,\mathcal A}\rightarrow(\mathbb R,\mathcal B)$ is called the \emph{experienced variable}, if
\begin{equation}\label{random-variable1}
{\bra{\xi_{\!\mathscr R}}}^{-1}(B)\in\bra{\mathcal A}
\end{equation}
for any Borel set $B\in\mathcal B$, i.e., a set  ${\bra{\xi_{\!\mathscr R}}}^{-1}(B)$ is a bra-event.
Equivalently speaking, the function ${\bra{\xi_{\!\mathscr R}}}={\bra{\xi_{\!\mathscr R}}}(\bra{\omega})$, which defined on the bra-set $\bra{\Omega}$ with values in $\mathbb R$, is called the \emph{experienced variable}, if
\begin{equation}\label{random-variable2}
\Big\{\bra{\omega}\colon {\bra{\xi_{\!\mathscr R}}}(\bra{\omega})<r\Big\} \in \bra{\mathcal A}
\end{equation}
for every choice of a real number $r\in\mathbb{R}$,
in other words, the set of elementary bra-incomes $\bra{\omega}$ such that ${\bra{\xi_{\!\mathscr R}}}(\bra{\omega})<r$ belongs to the bra-algebra $\bra{\mathcal A}$.

\texttt{\linebreak\indent Example \!\refstepcounter{ctrexa}\arabic{ctrexa}\,\label{exa-experienced-variable}\itshape\footnotesize (probability of ket-events as an experienced variable)\!.} Probabilities $p_x=\mathbf P(\ket{x})$ of ket-events  $\ket{x}\subseteq\ket{\Omega}$, $x\in\frak X_{\!\mathscr R}$ define on the believability space $\bra{\Omega}$ the function ${\bra{\xi_{\!\mathscr R}}}$ that takes on ech dual bra-event $\bra{x}\subseteq\bra{\Omega}$ the corresponding constant value
\begin{equation}\label{random-variable2}
{\bra{\xi_{\!\mathscr R}}}(\bra{\omega}) = p_x
\end{equation}
for all $\bra{\omega}\in\bra{x}$.
Since $\displaystyle\bra{\Omega}=\sum_{x\in\frak X _{\!\mathscr R}}\bra{x}$ then the function ${\bra{\xi_{\!\mathscr R}}}$ is defined on all the  bra-set $\bra{\Omega}$ and for any Borel $B\in\mathcal B$ the set ${\bra{\xi_{\!\mathscr R}}}^{-1}(B)$ is a bra-event:
\begin{equation}\label{random-variable2}
{\bra{\xi_{\!\mathscr R}}}^{-1}(B)=\sum_{{x\in\frak X _{\!\mathscr R}}\atop{p_x\in B}}\bra{x}\in\bra{\mathcal A},
\end{equation}
and the function  ${\bra{\xi_{\!\mathscr R}}}$ is the experienced variable by Definition \ref{def-experienced-variable}.

%\clearpage
\subsection{Random variable}%\\ \tiny $\sim$ \textsf{experienced-and-random-variables, 10 May 2017}}

\texttt{\linebreak\indent Definition \!\refstepcounter{ctrdef}\arabic{ctrdef}\,\label{def-random-variable}\itshape\footnotesize (random variable)\!.}
The function ${\ket{\xi_{\!\mathscr R}}}:\ket{\Omega,\mathcal A}\rightarrow(\mathbb R,\mathcal B)$ is called the \emph{random variable} if \begin{equation}\label{random-variable1}
{\ket{\xi_{\!\mathscr R}}}^{-1}(B)\in\ket{\mathcal A}
\end{equation}
for any Borel set $B\in\mathcal B$, i.e., a set  ${\ket{\xi_{\!\mathscr R}}}^{-1}(B)$ is a ket-event.
Equivalently speaking, the function  ${\ket{\xi_{\!\mathscr R}}}={\ket{\xi_{\!\mathscr R}}}(\ket{\omega})$, which is defined on the ket-set $\ket{\Omega}$ with values in $\mathbb R$, is called the \emph{random variable} if
\begin{equation}\label{random-variable2}
\Big\{\ket{\omega}\colon {\ket{\xi_{\!\mathscr R}}}(\ket{\omega})<r\Big\} \in \ket{\mathcal A}
\end{equation}
for every choice of a real number $r\in\mathbb{R}$,
in other words, the set of elementary ket-outcomes $\ket{\omega}$ such that ${\ket{\xi_{\!\mathscr R}}}(\ket{\omega})<r$ belongs to the ket-algebra $\ket{\mathcal A}$.

\texttt{\linebreak\indent Example \!\refstepcounter{ctrexa}\arabic{ctrexa}\,\label{exa-random-variable}\itshape\footnotesize (believability of bra-events as a random variable)\!.} The believability $b(X)=\mathbf B(\bra{\textsf{Ter}_{X/\!\!/\frak X_{\!\mathscr R}}})$ of terraced bra-events  $\bra{\textsf{Ter}_{X/\!\!/\frak X_{\!\mathscr R}}}\subseteq\bra{\Omega}$, $X\in\reflectbox{\bf S}^{\frak X _{\!\mathscr R}}$ define on the probability space $\ket{\Omega}$ the function ${\ket{\xi_{\!\mathscr R}}}$, which takes on every dual terraced ket-event $\ket{\textsf{ter}(X/\!\!/\frak X_{\!\mathscr R})}\subseteq\ket{\Omega}$ the corresponding constant value
\begin{equation}\label{random-variable2}
{\ket{\xi_{\!\mathscr R}}}(\ket{\omega}) = b(X)
\end{equation}
for all $\ket{\omega}\in\ket{\textsf{ter}(X/\!\!/\frak X_{\!\mathscr R})}$.
Since $\displaystyle\ket{\Omega}=\sum_{X\in\,\reflectbox{\scriptsize\bf S}^{\frak X _{\!\mathscr R}}}\ket{\textsf{ter}(X/\!\!/\frak X_{\!\mathscr R})}$ the the function ${\ket{\xi_{\!\mathscr R}}}$ is defined ob all the ket-set $\ket{\Omega}$ and for any Borel $B\in\mathcal B$ the set ${\ket{\xi_{\!\mathscr R}}}^{-1}(B)$ is a ket-event:
\begin{equation}\label{random-variable2}
{\ket{\xi_{\!\mathscr R}}}^{-1}(B)=\sum_{{X\in\,\reflectbox{\scriptsize\bf S}^{\frak X _{\!\mathscr R}}}\atop{b(X)\in B}}\ket{\textsf{ter}(X/\!\!/\frak X_{\!\mathscr R})}\in\ket{\mathcal A},
\end{equation}
and the function  ${\ket{\xi_{\!\mathscr R}}}$ is a random variable by Definition \ref{def-random-variable}.

\subsection{Experienced-random variable}%\\ \tiny $\sim$ \textsf{experienced-and-random-variables, 10 May 2017}}

\texttt{\linebreak\indent Definition \!\refstepcounter{ctrdef}\arabic{ctrdef}\,\label{def-experienced-random-variable}\itshape\footnotesize (experienced-random variable)\!.}
The function ${\braket{\xi_{\!\mathscr R}|\xi_{\!\mathscr R}}}:\braket{\Omega,\mathcal A|\Omega,\mathcal A}\rightarrow(\mathbb R,\mathcal B)$ is called the \emph{experienced-random variable} if \begin{equation}\label{random-variable1}
{\braket{\xi_{\!\mathscr R}|\xi_{\!\mathscr R}}}^{-1}(B)\in\braket{\mathcal A|\mathcal A}
\end{equation}
for any Borel set $B\in\mathcal B$, i.e. a set  ${\braket{\xi_{\!\mathscr R}|\xi_{\!\mathscr R}}}^{-1}(B)$ is a bra-ket-event.
Equivalently speaking, the function ${\braket{\xi_{\!\mathscr R}|\xi_{\!\mathscr R}}}={\braket{\xi_{\!\mathscr R}|\xi_{\!\mathscr R}}}(\braket{\omega^*|\omega})$ defined on the bra-ket-set $\braket{\Omega|\Omega}$ with values in $\mathbb R$ is called the \emph{experienced-random variable} if
\begin{equation}\label{random-variable2}
\Big\{\braket{\omega^*|\omega}\colon {\braket{\xi_{\!\mathscr R}|\xi_{\!\mathscr R}}}(\braket{\omega^*|\omega})<r\Big\} \in \braket{\mathcal A|\mathcal A}
\end{equation}
for every choice of a real number $r\in\mathbb{R}$,
in other words, the set of elementary bra-ket-incomes-outcomes $\braket{\omega^*|\omega}$ such that ${\braket{\xi_{\!\mathscr R}|\xi_{\!\mathscr R}}}(\braket{\omega^*|\omega})<r$ belongs to the bra-ket-algebra $\braket{\mathcal A|\mathcal A}$.

\texttt{\linebreak\indent Example \!\refstepcounter{ctrexa}\arabic{ctrexa}\,\label{exa-experienced-random-variable}\itshape\footnotesize (certainty of bra-ket-events as an experienced-random variable)\!.} Certainties $\bfphi_x(X/\!\!/\frak X_{\!\mathscr R})=\bfPhi(\braket{x|\textsf{ter}(X/\!\!/\frak X_{\!\mathscr R})})$ of bra-ket-events  $\braket{x|\textsf{ter}(X/\!\!/\frak X_{\!\mathscr R})}\subseteq\braket{\Omega|\Omega}$, $x\in\frak X_{\!\mathscr R}, X\in\reflectbox{\bf S}^{\frak X_{\!\mathscr R}}$ define on the certainty space $\braket{\Omega|\Omega}$ the function ${\braket{\xi_{\!\mathscr R}|\xi_{\!\mathscr R}}}$ which takes on each bra-ket-event $\braket{x|\textsf{ter}(X/\!\!/\frak X_{\!\mathscr R})}\subseteq\braket{\Omega|\Omega}$ the corresponding constant value
\begin{equation}\label{random-variable2}
{\braket{\xi_{\!\mathscr R}|\xi_{\!\mathscr R}}}(\braket{\omega^*|\omega}) = \bfphi_x(X/\!\!/\frak X_{\!\mathscr R})
\end{equation}
for all $\braket{\omega^*|\omega}\in\braket{x|\textsf{ter}(X/\!\!/\frak X_{\!\mathscr R})}$.
Since $\displaystyle\braket{\Omega|\Omega}=\sum_{x\in\frak X _{\!\mathscr R}}\sum_{X\in\,\reflectbox{\scriptsize\bf S}^{\frak X _{\!\mathscr R}}}\braket{x|\textsf{ter}(X/\!\!/\frak X_{\!\mathscr R})}$ then the function ${\braket{\xi_{\!\mathscr R}|\xi_{\!\mathscr R}}}$ is defined on all the bra-ket-set $\braket{\Omega|\Omega}$ and for any Borel $B\in\mathcal B$ set ${\braket{\xi_{\!\mathscr R}|\xi_{\!\mathscr R}}}^{-1}(B)$ is a  bra-ket-event:
\begin{equation}\label{random-variable2}
{\braket{\xi_{\!\mathscr R}|\xi_{\!\mathscr R}}}^{-1}(B)=\sum_{{(x,X)\in\frak X _{\!\mathscr R}\times\,\reflectbox{\scriptsize\bf S}^{\frak X _{\!\mathscr R}}}\atop{\bfphi_x(X/\!\!/\frak X _{\!\mathscr R})\in B}}\braket{x|\textsf{ter}(X/\!\!/\frak X_{\!\mathscr R})}\in\braket{\mathcal A|\mathcal A},
\end{equation}
and the function  ${\braket{\xi_{\!\mathscr R}|\xi_{\!\mathscr R}}}$ is an experienced-random variable by Definition \ref{def-experienced-random-variable}.

\texttt{\linebreak\indent Definition \!\refstepcounter{ctrdef}\arabic{ctrdef}\,\label{def-dfs}\itshape\footnotesize (functions of distributions of believabilities, probabilities and certainties)\!.} The functions
\begin{equation}\label{dfs}
\begin{split}
F_{\bra{\xi_{\!\mathscr R}}}(r) &= \mathbf B(\{\bra{\omega}\colon {\bra{\xi_{\!\mathscr R}}}(\bra{\omega})<r\}) = \mathbf B({\bra{\xi_{\!\mathscr R}}}<r),\\
F_{\ket{\xi_{\!\mathscr R}}}(r) &= \mathbf P(\{\ket{\omega}\colon {\ket{\xi_{\!\mathscr R}}}(\ket{\omega})<r\}) = \mathbf P({\ket{\xi_{\!\mathscr R}}}<r),\\
F_{\braket{\xi_{\!\mathscr R}|\xi_{\!\mathscr R}}}(r) &= \bfPhi(\{\braket{\omega^*|\omega}\colon {\braket{\xi_{\!\mathscr R}|\xi_{\!\mathscr R}}}(\braket{\omega^*|\omega})<r\}) = \bfPhi({\braket{\xi_{\!\mathscr R}|\xi_{\!\mathscr R}}}<r),
\end{split}
\end{equation}
where $-\infty$ and $+\infty$ are allowed as values $r$, are called
\emph{the function of  believability distribution of the experienced variable} ${\bra{\xi_{\!\mathscr R}}}$,
\emph{the function of  probability distribution of the random variable} ${\ket{\xi_{\!\mathscr R}}}$,
 and \emph{the function of  certainty distribution of the experienced-random variable} ${\braket{\xi_{\!\mathscr R}|\xi_{\!\mathscr R}}}$ correspondingly.

\section{Dual inducing the nonadditive functions of a set by believability and probability\label{versus}}%\\ \tiny $\sim$ \textsf{0002/dual-problems}}

In the theory of experience and chance for each co$\sim$event $\mathscr R\subseteq\braket{\Omega|\Omega}$
the \emph{believability} $\mathbf B$ defined on sigma-algebra of believability space $\bra{\Omega,\mathcal A,\mathbf B}$ $\mathscr R$-induces on the probability space, and the \emph{probability} $\mathbf P$ defined on  sigma-algebra of probability space $\ket{\Omega,\mathcal A,\mathbf P}$  $\mathscr R$-induces on the believability space the functions of a set which do not possess a property of additivity on these spaces.

Let us consider this fact in more detail, since it is, in my opinion, for a long time misleads the apologists of fuzzy mathematics \cite{Zadeh1965,Shafer1976,DuboisPrade2001,Zimmermann2001,Zadeh2014} and forces them in their articles to make mandatory statements that those set functions that they intend to deal with (possibilities, beliefs, etc.) are not absolutely a \emph{probability}, so as they do not have the property of additivity, and they are \emph{not related to a probability}.
The origins of these misconceptions are outlined in my works \cite{Vorobyev2009em, Vorobyev2009a}.
Now my explanations of this aberration are based entirely on the axiomatics of the theory of experience and chance and consist in the following. Those set functions that are of interest in fuzzy math are always so or otherwise $\mathcal R$-induced by the probability the set functions on the believability space, which, naturally, do not possess an additivity of this space, but in the theory of experience and chance are mutually unambiguously associated with the additive on this space the set function, which I also call a \emph{believability}. In the theory of experience and chance, the dual assertion is also true: the believability measure, additive on the believability space, in its turn $\mathscr R$ -induces on the probability space, the non-additive functions of a set that are one-to-one related to the probability. So, consider the relationships between $\mathscr R$-induced nonadditive functions of a set on the one hand and a \emph{believability} and a \emph{probability} on the other.

The co$\sim$event $\mathscr R\subseteq\braket{\Omega|\Omega}$ by probability $\mathbf P$ defined on $\ket{\Omega}$,  and by believability $\mathbf B$ defined on $\bra{\Omega}$ induces:
\begin{itemize}
\item
on $\bra{\Omega,\mathcal A,\mathbf B}$ the nonadditive set function  $\mathbf P'$ defining its values on each bra-event $\bra{x}\subseteq\bra{\Omega}, x\in\frak X_{\!\mathscr R}$, dual to the ket-event $\ket{x}\subseteq\ket{\Omega}$, and on each terraced bra-event $\bra{\textsf{Ter}_{X/\!\!/\frak X_{\!\mathscr R}}}\subseteq\bra{\Omega}, X\in\reflectbox{\bf S}^{\frak X_{\!\mathscr R}}$, dual to the terraced ket-event $\ket{\textsf{ter}(X/\!\!/\frak X_{\!\mathscr R})}\subseteq\ket{\Omega}$, by the formulas:
\begin{equation}\label{P'}
\begin{split}
\mathbf P'(\bra{x})&=\mathbf P(\ket{x}),\\
\mathbf P'(\bra{\textsf{Ter}_{X/\!\!/\frak X_{\!\mathscr R}}})&=\mathbf P(\ket{\textsf{ter}(X/\!\!/\frak X_{\!\mathscr R})});
\end{split}
\end{equation}
\item
on $\ket{\Omega,\mathcal A,\mathbf P}$ the nonadditive set function $\mathbf B'$, defining its values on each ket-event $\ket{x}\subseteq\ket{\Omega}, x\in\frak X_{\!\mathscr R}$, dual to the bra-event $\bra{x}\subseteq\bra{\Omega}$, and on each terraced ket-event $\ket{\textsf{ter}(X/\!\!/\frak X_{\!\mathscr R})}\subseteq\ket{\Omega}, X\in\reflectbox{\bf S}^{\frak X_{\!\mathscr R}}$, dual to the terraced bra-event $\bra{\textsf{Ter}_{X/\!\!/\frak X_{\!\mathscr R}}}\subseteq\bra{\Omega}$, by formulas:
\begin{equation}\label{B'}
\begin{split}
\mathbf B'(\ket{x})&=\mathbf B(\bra{x}),\\
\mathbf B'(\ket{\textsf{ter}(X/\!\!/\frak X_{\!\mathscr R})})&=\mathbf B(\bra{\textsf{Ter}_{X/\!\!/\frak X_{\!\mathscr R}}}).
\end{split}
\end{equation}
\end{itemize}

\texttt{\linebreak\indent Property \!\refstepcounter{ctrATT}\arabic{ctrATT}\,\label{ATT-non-additivity}\itshape\footnotesize (non-additivity of induced set functions)\!.} \emph{The induced set functions $\mathbf P'$ and $\mathbf B'$ are not additive on $\bra{\Omega,\mathcal A,\mathbf B}$ and $\ket{\Omega,\mathcal A,\mathbf B}$ correspondingly.}

\texttt{Proof.} Since the probability $\mathbf P$ is additive on the ket-space $\ket{\Omega,\mathcal A,\mathbf B}$, and
$$
\ket{x}=\sum_{x\in X\in\,\reflectbox{\scriptsize\bf S}^{\frak X_{\!\mathscr R}}} \ket{\textsf{ter}(X/\!\!/\frak X_{\!\mathscr R})},
$$
then for $x\in\frak X_{\!\mathscr R}$
$$
\mathbf P(\ket{x})=\sum_{x\in X\in\,\reflectbox{\scriptsize\bf S}^{\frak X_{\!\mathscr R}}} \mathbf P(\ket{\textsf{Ter}_{X/\!\!/\frak X_{\!\mathscr R}}}).
$$
From this and (\ref{P'}) we get that
$$
\mathbf P'(\bra{x})=\sum_{x\in X\in\,\reflectbox{\scriptsize\bf S}^{\frak X_{\!\mathscr R}}} \mathbf P'(\bra{\textsf{Ter}_{X/\!\!/\frak X_{\!\mathscr R}}}),
$$
but since for $X\in\reflectbox{\bf S}^{\frak X_{\!\mathscr R}}$
$$
\bra{\textsf{Ter}_{X/\!\!/\frak X_{\!\mathscr R}}}=\sum_{x\in X\in\,\reflectbox{\scriptsize\bf S}^{\frak X_{\!\mathscr R}}} \bra{x},
$$
then, generally speaking,
$$
\bra{x}\ne\sum_{x\in X\in\,\reflectbox{\scriptsize\bf S}^{\frak X_{\!\mathscr R}}} \bra{\textsf{Ter}_{X/\!\!/\frak X_{\!\mathscr R}}},
$$
which proves the non-additivity of the induced set function $\mathbf P'$ on $\bra{\Omega,\mathcal A,\mathbf B}$.
Similarly, since the believability $\mathbf B$ is additive on the bra-space $\bra{\Omega,\mathcal A,\mathbf B}$, and
$$
\bra{\textsf{Ter}{X/\!\!/\frak X_{\!\mathscr R}}}=\sum_{x\in X\in\,\reflectbox{\scriptsize\bf S}^{\frak X_{\!\mathscr R}}} \bra{x},
$$
then for $X\in\reflectbox{\bf S}^{\frak X_{\!\mathscr R}}$
$$
\mathbf B(\bra{\textsf{Ter}_{X/\!\!/\frak X_{\!\mathscr R}}})=\sum_{x\in X\in\,\reflectbox{\scriptsize\bf S}^{\frak X_{\!\mathscr R}}} \mathbf B(\bra{x}).
$$
From this and (\ref{B'}) we get that
$$
\mathbf B'(\ket{\textsf{ter}(X/\!\!/\frak X_{\!\mathscr R})})=\sum_{x\in X\in\,\reflectbox{\scriptsize\bf S}^{\frak X_{\!\mathscr R}}} \mathbf B'(\ket{x}),
$$
but since for $x\in\frak X_{\!\mathscr R}$
$$
\ket{x}=\sum_{x\in X\in\,\reflectbox{\scriptsize\bf S}^{\frak X_{\!\mathscr R}}} \ket{\textsf{ter}(X/\!\!/\frak X_{\!\mathscr R})},
$$
then, generally speaking,
$$
\ket{\textsf{ter}(X/\!\!/\frak X_{\!\mathscr R})}\ne\sum_{x\in X\in\,\reflectbox{\scriptsize\bf S}^{\frak X_{\!\mathscr R}}} \ket{x},
$$
which proves the non-additivity of induced set function $\mathbf B'$ on $\ket{\Omega,\mathcal A,\mathbf B}$.

\clearpage
\section{Examples of the use of certainty theory}

We will mention only two examples of the use of the new theory of experience and chance, one of which (``student delights'') is discussed in this article and shows for the time being only a curious connection between the two dualities: in the theory of experience and chance and in the theory of optimization; and the second (``the bet on a bald'') is discussed in detail in my other work \cite{Vorobyev2016famems5} and is devoted to the correct mathematical description of experienced-random experiment, which, although carried out at the macro level, but in which the observer clearly affects the outcome of the observation accurately just as in physics at the quantum level.

\subsection{Problem of ``student delicacies''}

The student decides which purchase to make in the bakery for an after-dinner delicacy.
There is the set $\frak X_{\!\mathscr R}$ of delicacies $x\in\frak X_{\!\mathscr R}$.
The delicacies contain healthy ingredients $X\in\reflectbox{\bf S}^{\frak X_{\!\mathscr R}}$ forming the set $\reflectbox{\bf S}^{\frak X_{\!\mathscr R}}$.
``The delicacy  $x\in\frak X_{\!\mathscr R}$ the student buys, i.e., the ket-event $\ket{x}\subseteq\ket{\Omega}$  happens'' with probability $p_x$.
Taking care of his health, the student decided that his believability
``in the benefits of ingredients (in terraced bra-events)'' $\bra{\textsf{Ter}_{X/\!\!/\frak X_{\!\mathscr R}}}\subseteq\bra{\Omega},\,X\in\reflectbox{\bf S}^{\frak X_{\!\mathscr R}}$ should be at least $b(X/\!\!/\frak X_{\!\mathscr R})$:
\begin{equation}\label{min-subject-to}
\begin{split}
&\sum_{x\in X} b_x \geqslant b(X/\!\!/\frak X),
\end{split}
\end{equation}
where $b_x$ is the believability ``in the benefits of delicacies (in bra-events)'' $\bra{x}\subseteq\bra{\Omega},\,x\in\frak X_{\!\mathscr R}$ for her/his health.

The problem of ``student delicacies'' can be formulated as follows:
\begin{equation}\label{min}
\begin{split}
\min_{b_{x},\,x\in\frak X_{\!\mathscr R}}\hspace{10pt}& \sum_{x\in\frak X_{\!\mathscr R}} b_{x}p_{x}\\
\mbox{subject to}\hspace{10pt}&\\[-16pt]
&\sum_{x\in X} b_x \geqslant b(X/\!\!/\frak X),\,X\in\,\reflectbox{\bf S}^{\frak X_{\!\mathscr R}},\\
&\sum_{x\in\frak X_{\!\mathscr R}} b_{x}=1,\,b_{x}\geqslant 0,\,x\in \frak X_{\!\mathscr R}
\end{split}
\end{equation}
--- \textbf{she/he is looking for a believability distribution $\{b_x\colon\in\frak X_{\!\mathscr R}\}$ on which the mean-believable probability $\bfEr_{b_x}(p_x)=\sum_{x\in\frak X_{\!\mathscr R}} b_{x}p_{x}$ of ``purchases of delicacies (the ket-events)'' $\ket{x}\subseteq\ket{\Omega},\,x\in\frak X_{\!\mathscr R}$ takes a minimal value under the constraints (\ref{min-subject-to}) made.}
Here $b_{x}$ is the believability ``in the benefits of delicacies (in bra-events)'' $\bra{x}\subseteq\bra{\Omega},\,x\in\frak X_{\!\mathscr R}$ for her/his health;
$p_{x}$ is the probability of ``purchases of delicacies (the ket-events)'' $\ket{x}\subseteq\ket{\Omega},\,x\in\frak X_{\!\mathscr R}$;
$b(X/\!\!/\frak X_{\!\mathscr R})$ is the believability ``in the benefits of delicacies (in terraced bra-events)'' $\bra{\textsf{Ter}_{X/\!\!/\frak X_{\!\mathscr R}}}\subseteq\bra{\Omega},\,X\in\,\reflectbox{\bf S}^{\frak X_{\!\mathscr R}}$.

\begin{table}[h!]
\tiny
\centering
\begin{tabular}{|c|c|c|c|c|c|}
  \hline
                             &believability in the benefits         &ingredient $X_1$&ingredient $X_2$&ingredient $X_3$&probability of purchases\\
                             & of delicacies $x\in\frak X_{\!\mathscr R}$&&&&of delicacies $x\in\frak X_{\!\mathscr R}$                                        \\ \hline
probability of purchases          &                     & & & &                                                         \\
a set of delicacies $X\in\reflectbox{\bf S}^{\frak X_{\!\mathscr R}}$&                     & $p(X_1)=0.2$& $p(X_2)=0.3$&$p(X_3)=0.5$&                      \\ \hline
delicacy $x_1$              &$b_{x_1}=0.6$        & 0.12 $_\bullet$& 0.18 $_\bullet$& 0.30 $_\circ$& $p_{x_1}=0.5$                                  \\ \hline
delicacy $x_2$              &$b_{x_2}=0.4$        & 0.08 $_\circ$& 0.12 $_\bullet$& 0.20 $_\bullet$& $p_{x_2}=0.8$                                  \\ \hline
believability in the benefit         &                     & & & &                                                         \\
of ingredients $X\in\reflectbox{\bf S}^{\frak X_{\!\mathscr R}}$    &                     & $b(X_1)=0.6$& $b(X_2)=1.0$& $b(X_3)=0.4$&                     \\ \hline
\end{tabular}
\caption{Data for the problem ``student delicacies'' with 2 delicacies forming the doublet $\frak X_{\!\mathscr R}=\{x_1,x_2\}$, and 3 ingredients forming the triplet $\protect\reflectbox{\bf S}^{\frak X_{\!\mathscr R}}=\{X_1,X_2,X_3\}=\{\{x_1\},\{x_1,x_2\},\{x_3\}\}$. Table cells marked with black circles correspond to the co$\sim$event $\mathscr R=\braket{x_1|\textsf{ter}(X_1/\!\!/\frak X_{\!\mathscr R})}+\braket{x_1|\textsf{ter}(X_2/\!\!/\frak X_{\!\mathscr R})}+\braket{x_2|\textsf{ter}(X_2/\!\!/\frak X_{\!\mathscr R})}+\braket{x_2|\textsf{ter}(X_3/\!\!/\frak X_{\!\mathscr R})}\subseteq\braket{\Omega|\Omega}$.\label{table-LP}}
\end{table}

We now accept the dualistic point of view of the student, from which she/he places a restriction on the probability of ``purchases of delicacies (the ket-events)'' $\ket{x}\subseteq\ket{\Omega},\,x\in\frak X_{\!\mathscr R}$:
\begin{equation}\label{max-subject-to}
\begin{split}
\sum_{x\in X\in\,\reflectbox{\scriptsize\bf S}^{\frak X_{\!\mathscr R}}} p(X/\!\!/\frak X_{\!\mathscr R}) \leqslant p_x,
\end{split}
\end{equation}
where $p(X/\!\!/\frak X_{\!\mathscr R})$ is the probability of ``purchaces of ingredients (the terraced ket-events)'' $\ket{\textsf{ter}(X/\!\!/\frak X_{\!\mathscr R})}\subseteq\ket{\Omega},\,X\in\,\reflectbox{\bf S}^{\frak X_{\!\mathscr R}}$.

The dual problem solved by the student is as follows:
\begin{equation}\label{max}
\begin{split}
\max_{p(X),\,X\in\,\reflectbox{\scriptsize\bf S}^{\frak X_{\!\mathscr R}}}\hspace{10pt}& \sum_{X\in\,\reflectbox{\scriptsize\bf S}^{\frak X_{\!\mathscr R}}} b(X/\!\!/\frak X_{\!\mathscr R})p(X/\!\!/\frak X_{\!\mathscr R})\\
\mbox{subject to}\hspace{10pt}&\\[-16pt]
&\sum_{x\in X\in\,\reflectbox{\scriptsize\bf S}^{\frak X_{\!\mathscr R}}}p(X/\!\!/\frak X_{\!\mathscr R}) \leqslant p_{x},\\
&\sum_{X\in\,\reflectbox{\scriptsize\bf S}^{\frak X_{\!\mathscr R}}}p(X/\!\!/\frak X_{\!\mathscr R})=1,\,p(X/\!\!/\frak X_{\!\mathscr R}) \geqslant 0,\, X\in\,\reflectbox{\scriptsize\bf S}^{\frak X_{\!\mathscr R}}
\end{split}
\end{equation}
--- \textbf{the student looks for a probability distribution $\{p(X/\!\!/\frak X_{\!\mathscr R})\colon X\in\reflectbox{\bf S}^{\frak X_{\!\mathscr R}}\}$ on which the mean-prabable believability $\bfEl_{p(X)}(b(X))=\sum_{X\in\,\reflectbox{\scriptsize\bf S}^{\frak X_{\!\mathscr R}}} b(X/\!\!/\frak X_{\!\mathscr R})p(X/\!\!/\frak X_{\!\mathscr R})$ ``in the benefits for her/his health of ingredients (in terraced bra-events)'' $\bra{\textsf{Ter}_{X/\!\!/\frak X_{\!\mathscr R}}}\subseteq\bra{\Omega},\,X\in\reflectbox{\bf S}^{\frak X_{\!\mathscr R}}$ takes a maximal value under the constraints (\ref{max-subject-to}) made.}
Here $p(X/\!\!/\frak X_{\!\mathscr R})$ is the probability of ``purchases of ingredients (the terraced ket-events)'' $\ket{\textsf{ter}(X/\!\!/\frak X_{\!\mathscr R})}\subseteq\ket{\Omega},\,X\in\,\reflectbox{\bf S}^{\frak X_{\!\mathscr R}}$;
$p_{x}$ is the probability of ``purchases of delicacies (the ket-events)'' $\ket{x}\subseteq\ket{\Omega},\,x\in\frak X_{\!\mathscr R}$;
$b(X/\!\!/\frak X_{\!\mathscr R})$ is the believability ``in the benefits of ingredients  (in terraced bra-events)'' $\bra{\textsf{Ter}_{X/\!\!/\frak X_{\!\mathscr R}}}\subseteq\bra{\Omega},\,X\in\,\reflectbox{\scriptsize\bf S}^{\frak X_{\!\mathscr R}}$.

In the matrix form, the direct problem can be expressed as: \emph{``To minimize $p^T\breve{b}$ under the condition $A\breve{b} \geqslant b$, $\breve{b} \geqslant 0$, $\breve{\mathbf I}^T\breve{b}=1$''}; with the corresponding dual problem: \emph{``To minimize $\breve{b}^Tp$ under the condition $A^Tp \leqslant \breve{p}$, $p \geqslant 0$, $\mathbf I^Tp=1$''};
where $\breve{b}=\{b_x\colon x\in\frak X_{\!\mathscr R}\}$,
$\breve{p}=\{p_x\colon x\in\frak X_{\!\mathscr R}\}$,
$b=\left\{b(X/\!\!/\frak X_{\!\mathscr R})\colon X\in\reflectbox{\bf S}^{\frak X_{\!\mathscr R}}\right\}$,
$p=\left\{p(X/\!\!/\frak X_{\!\mathscr R})\colon X\in\reflectbox{\bf S}^{\frak X_{\!\mathscr R}}\right\}$,
$\breve{\mathbf I}=\{1\colon x\in\frak X_{\!\mathscr R}\}$,
$\mathbf I=\left\{1\colon X\in\reflectbox{\bf S}^{\frak X_{\!\mathscr R}}\right\}$ are set-columns, $\breve{b}^T,\breve{p}^T,b^T,p^T,\mathbf I^T$ are set-rows correspondingly, and $A=\left\{\mathbf 1_X(x)\colon x\in\frak X_{\!\mathscr R}, X\in\reflectbox{\bf S}^{\frak X_{\!\mathscr R}}\right\}$ is the set-matrix.

\section{Instead of discussing}

Before the finish I have to slow down on three sharp corners.

On the first one, we need to stop and carefully study the main innovation of this work, Axiom \ref{Ax-co-event} on page \pageref{Ax-co-event}, which expands the ``silent'' Kolmogorov axiom of an event, so that this axiom together with its dual reflection allowed a new theory to jointly explore both the future randomness of observations, and the past experience of observers.

On the second, it is impossible to rush past the very curious temporal bra-ket-duality of statements from the property \ref{ATT-bra-ket-events2} on the page \pageref{ATT-bra-ket-events2}, which states that
\begin{itemize}
\item
from the Kolmogorov theory of probabilities:
If there is some ket-event, then with it all ket-events occur, in which it is contained as a ket-subset. In other words, all the ket-events which follow from it, i.e. which can serve as \emph{its consequences in the future}.
\item
from the dual theory of believabilities:
If some bra-event is experienced then with it are experienced all the bra-events, which it contains as bra-subsets. In other words, all bra-events from which it follows, i.e.,
which could serve as \emph{its causes in the past}.
\end{itemize}
This remarkable property of temporal duality ket-events and bra-events clearly shows the similarity and difference between the \emph{future chance} and the \emph{past experience}, which for the first time are jointly mathematically correctly studied in the theory experience and chance postulated in this article (see also \cite{Vorobyev2016famems3}).

And finally, on the third one, it is worthwhile once again to linger on explaining the new theory (see Property \ref{ATT-non-additivity} on the page \pageref{ATT-non-additivity}) for quite a long time confusing the apologists of fuzzy mathematics \cite{Zadeh1965,Shafer1976,DuboisPrade2001,Zimmermann2001,Zadeh2014}
on the non-additivity of the set functions of interest, the origins of which are considered in my works earlier \cite[2009]{Vorobyev2009em,Vorobyev2009a}.

{\footnotesize
\bibliography{vorobyev}

\begin{thebibliography}{10}

\bibitem{Vorobyev2016famems1}
O.~Yu. Vorobyev.
\newblock An element-set labelling a {C}artesian product by measurable binary
  relations which leads to postulates of the theory of experience and chance as
  a theory of co$\sim$events.
\newblock {\em \emph{In.} Proc. of the XV Intern. FAMEMS Conf. on Financial and
  Actuarial Mathematics and Eventology of Multivariate Statistics \& the
  Workshop on Hilbert's Sixth Problem; Krasnoyarsk, SFU (Oleg Vorobyev ed.)},
  pages 9--24, ISBN 978--5--9903358--6--8,
  \mbox{\tiny\url{https://www.academia.edu/34390291}}, 2016.

\bibitem{Dirac1939}
P.A.M. Dirac.
\newblock A new notation for quantum mechanics.
\newblock {\em Proceedings of the Cambridge Philosophical Society},
  (35):416--418, 1939.

\bibitem{Dirac1964}
P.A.M. Dirac.
\newblock {\em The principles of quantum mechanics}.
\newblock Yeshiva University, New York, 1964.

\bibitem{Vorobyev2007}
O.~Yu. Vorobyev.
\newblock {\em Eventology}.
\newblock Siberian Federal University, Krasnoyarsk, Russia, 435p.,
  \mbox{\tiny\url{https://www.academia.edu/179393/}}, 2007.

\bibitem{Wittgenstein1921}
L.~Wittgenstein.
\newblock Logisch-philosophische abhandlung.
\newblock {\em Ostwalds Annalen der Naturphilosophie}, 14:185--262, 1921.

\bibitem{Robbins1944/45}
H.~E. Robbins.
\newblock On the measure of a random set.
\newblock {\em Ann. Math. Statist.}, 15/16:70--74/342--347, 1944/45.

\bibitem{Fubini1907}
G.~Fubini.
\newblock Sugli integrali multipli.
\newblock {\em Rom. Acc. L. Rend.}, (5) 16(1):608--614, 1907.

\bibitem{Kolmogorov1933}
A.~N. Kolmogorov.
\newblock {\em Grundbegriffe der Wahrscheinlichkeitrechnung}.
\newblock Ergebnisse der Mathematik, Berlin, 1933.

\bibitem{Zadeh1965}
L.~A. Zadeh.
\newblock Fuzzy sets.
\newblock {\em Information and Control}, 8 (3):338--353, 1965.

\bibitem{Shafer1976}
G.~Shafer.
\newblock {\em A Mathematical Theory of Evidence}.
\newblock Princeton University Press, Princeton, NJ, 1976.

\bibitem{DuboisPrade2001}
D.~Dubois and H.~Prade.
\newblock Possibility theory, probability theory and multiple-valued logics: A
  clarification.
\newblock {\em Annals of Mathematics and Artificial Intelligence}, 32:35--66,
  2001.

\bibitem{Zimmermann2001}
H.-J. Zimmermann.
\newblock {\em Fuzzy Set Theory and Its Applications}.
\newblock Kluwer Academic Publishers, Dordrecht, the Netherlands, 2001.

\bibitem{Zadeh2014}
L.~A. Zadeh.
\newblock A note on similarity-based definitions of possibility and
  probability.
\newblock {\em Information Sciences}, 267:334--336, 2014.

\bibitem{Vorobyev2009em}
O.~Yu. Vorobyev.
\newblock Eventology versus contemporary theories of uncertainty.
\newblock {\em \emph{In.} Proc. of the XII Intern. {EM} conference on
  eventological mathematics and related fields, Krasnoyarsk: SFU (Oleg Vorobyev
  ed.)}, pages 13--27,
  \mbox{\tiny\url{https://mpra.ub.uni--muenchen.de/13961/1/MPRA\_paper\_13961.pdf}},
  2009.

\bibitem{Vorobyev2009a}
O.~Yu. Vorobyev.
\newblock Contemporary uncertainty theories: An eventological view.
\newblock {\em \emph{In.} Proc. of the VIII Intern. FAM Conf. on Financial and
  Actuarial Mathematics and Related Fields, Krasnoyarsk: SFU (Oleg Vorobyev
  ed.)}, 1:26--34, 2009.

\bibitem{Vorobyev2016famems5}
O.~Yu. Vorobyev.
\newblock The bet on a bald.
\newblock {\em \emph{In.} Proc. of the XV Intern. FAMEMS Conf. on Financial and
  Actuarial Mathematics and Eventology of Multivariate Statistics \& the
  Workshop on Hilbert's Sixth Problem; Krasnoyarsk, SFU (Oleg Vorobyev ed.)},
  pages 98--101, ISBN 978--5--9903358--6--8,
  \mbox{\tiny\url{https://www.academia.edu/33154290}}, 2016.

\bibitem{Vorobyev2016famems3}
O.~Yu. Vorobyev.
\newblock Theory of dual co$\sim$event means.
\newblock {\em \emph{In.} Proc. of the XV Intern. FAMEMS Conf. on Financial and
  Actuarial Mathematics and Eventology of Multivariate Statistics \& the
  Workshop on Hilbert's Sixth Problem; Krasnoyarsk, SFU (Oleg Vorobyev ed.)},
  pages 44--93, ISBN 978--5--9903358--6--8,
  \mbox{\tiny\url{https://www.academia.edu/34357251}}, 2016.

\end{thebibliography}
}

\end{document}